\documentclass[12pt,oneside]{article}
\usepackage{epsfig}
\usepackage{enumerate}
\usepackage{latexsym}
\usepackage{amsmath}
\usepackage{amssymb}
\usepackage[english]{babel}

\newtheorem{thm}{Theorem}[subsection]

\newtheorem{cor}[thm]{Corollary}
\newtheorem{lem}[thm]{Lemma}
\newtheorem{fact}[thm]{Fact}
\newtheorem{prop}[thm]{Proposition}
\newtheorem{defn}[thm]{Definition}

\newtheorem{expl}[thm]{Example}

\DeclareMathOperator*{\bigboxplus}{\text{\Huge $_{\boxplus}$}}
\DeclareMathOperator*{\Card}{\mathrm{Card}}

\newcommand{\Real}{\mathbb{R}}

\numberwithin{equation}{section}

\newfont{\sfl}{cmssi12}

\begin{document}

\title{On Some Idempotent and Non-Associative Convex Structure.}
 \date{ September 2013}

\author{\thanks {University of Perpignan, 52 avenue Villeneuve,
66000 Perpignan, France.}  Walter Briec }

\maketitle
\begin{abstract}
$\mathbb{B}$-convexity was defined in ~\cite{bh}
 as a suitable  Kuratowski-Painlev\'e upper limit of linear convexities
over a finite dimensional Euclidean vector space. Excepted in the
special case where convex sets are subsets of $\Real_+^n$,
$\mathbb B$-convexity was not defined with respect to a given
explicit algebraic structure. This is done in that paper, which
proposes an extension of $\mathbb B$-convexity to the whole
Euclidean vector space. An unital idempotent and non-associative
magma is defined over the real set and an extended $n$-ary
operation is introduced. Along this line, the existence of the
Kuratowski-Painlev\'e limit of the convex hull of two points over
$\mathbb R^n$ is shown and an explicit extension of $\mathbb
B$-convexity is proposed.
\end {abstract}

{\bf AMS:} 06D50, 32F17\\


\section{Introduction}\label{SECMXVSP}

 $\mathbb B$-convexity was defined in \cite{bh}. One can say, loosely
speaking, that this $\mathbb B$-convexity is obtained from the usual
linear convexity making the formal transformation $+\rightarrow
\max$.  By definition, a $\mathbb B$-convex subset of $\Real_+^n$ is
a connected upper semilattice. $\mathbb B$-convex functions were
analyzed in \cite{adil}. Hanh-Banach like separation properties
\cite{bh3} as well as fixed point results \cite{bh2} -see also
\cite{pat}- have been established. The standard form of $\mathbb
B$-convexity is defined on the nonegative Euclidean orthant
$\Real_+^n$ and is linked to Max-Plus algebra via a suitable
homeomorphism. In finite dimensional space, Max-Plus convexity and
$\mathbb B$-convexity are isomorphic topological Maslov's
semi-modules \cite{ms} and, consequently, a proposition that is true
in the framework of $\mathbb B$-convexity holds, with obvious
lexical modifications, in Max-Plus convexity. Though $\mathbb
B$-convexity was initially defined over $\Real^n$ as a
Kuratowski-Painlev\'e upper limit of linear convexities, it was not
described in term of an explicit algebraic structure, excepted in
the case where convex sets are subsets of $\Real_+^n$. More
recently, $\mathbb B^{-1}$-convex sets were introduced in
\cite{adilYe} and \cite{adilYeTi}.

This paper introduces a suitable algebraic structure  extending
$\mathbb B$-convexity to the whole Euclidean vector space. However,
there do not exist non trivial algebraic structures being both
idempotent, associative,
 and  having inverse elements. Therefore, a special class of idempotent magma is
 considered in which associativity  is relaxed to
preserve symmetry and idempotence. This binary operation is based
upon absolute value function. An $n$-ary extension of this
algebraic structure is proposed and related to the pointwise limit
of a generalized H\"older sum. Some algebraic properties are
established and an extended definition of $\mathbb B$-convexity is
then proposed, including as a special case that one proposed in
\cite{bh}. It is shown that such a notion of convexity is
equivalently characterized from the Kuratowski-Painlev\'e limit of
the generalized convex hull of two points defined in \cite{bh}.

The paper unfolds as follows. Section 2 focusses on a special class
of symmetrical idempotent magmas. An $n$-ary extension of this
operation is proposed and is related to the limit  of a generalized
H\"older sum. In section 3, it is established that such an algebraic
structure yields a very simple extension of $\mathbb B$-convexity
over $\mathbb R^n$.

\section{Pointwise Limit of a Generalized Sum and Algebraic Structure }\label{intralg}

For all $p\in \mathbb{N}$, we consider a bijection $\varphi_p
:\mathbb{R} \longrightarrow \mathbb{R}$ defined by:
$\varphi_p:\lambda\longrightarrow \lambda^{2p+1}$
 and $\Phi_p (x_1,\cdots ,x_n)=(\varphi_p(x_1),\cdots ,\varphi_p(x_n))$. We can induce a field
structure on $\mathbb{R}$ for which $\varphi_p$ becomes a field
isomorphism. Given this change of notation via $\varphi_p$ and
$\Phi_p$ we can define a $\Real$-vector space structure on
$\Real^n$ by: $\lambda\stackrel{\varphi_p}{.}
x=\Phi_p^{-1}(\varphi_p(\lambda).\Phi_p(x))=\lambda.x$ and $x_1
\stackrel{\varphi_p}{+}
x_2=\Phi_p^{-1}\big(\Phi_p(x_1)+\Phi_p(x_2)\big)$; we call these
two operations the indexed scalar product and the indexed sum
(indexed by $\varphi$ of course). For all natural number $n\geq
1$, let us denote $[n]=\{1,...,n\}$. Suppose that
$A=\{x_1,\cdots,x_m\}$. The $\varphi_p$-sum - denoted
$\stackrel{\varphi_p}{\sum}$ (of $(x_1,\cdots,x_m)\in
\mathbb{R}^n$)- is defined by
 $\stackrel{\varphi_p}{\sum}_{i\in [m]}x_i=
\Phi_p^{-1}\Big(\sum_{i\in [m]}\Phi_p (x_i)\Big)$. In the remainder
of the paper we denote for all $x,y\in \Real^n$
$x\stackrel{p}{+}y=x\stackrel{\varphi_p}{+}y$.

\subsection{On some Idempotent, Symmetrical and Non-associative Algebraic Structure}

If $x,y\in \Real_+$ then one has $\lim_{p\longrightarrow
+\infty}x\stackrel{p}{+}y=\max\{|x|,|y|\}$. In the case where
 $x$ and $y$ belong to the whole real set, it is easy to establish the
 following property.

\begin{fact} \label{two}For all $x,y\in \Real$ we have:
\begin{equation*}\label{base}\lim_{p\longrightarrow
+\infty}x\stackrel{p}{+}y=\lim_{p\longrightarrow
+\infty}\big(x^{2p+1}+y^{2p+1}\big)^{\frac{1}{2p+1}}=
\left\{\begin{matrix}x\ &\hbox{ if } &|x|&>&|y|\\
\frac{1}{2}(x+y)&\hbox{ if }&|x|&=&|y|\\
y& \hbox{ if }& |x|&<&|y|.\end{matrix}\right.\end{equation*}
\end{fact}
{\bf Proof:} Suppose that $x=-y$. Then for all $p\in \mathbb N$ we
have
$\big(x^{2p+1}+y^{2p+1}\big)^{\frac{1}{2p+1}}=\big(x^{2p+1}-x^{2p+1}\big)^{\frac{1}{2p+1}}=0$,
which proves  this case. If $x=y$, then $\lim_{p\longrightarrow
+\infty}\big(x^{2p+1}+y^{2p+1}\big)^{\frac{1}{2p+1}}=x\lim_{p\longrightarrow
+\infty}2^{\frac{1}{2p+1}}=x$. To end the proof, suppose, for
example, that $|x|>|y|$. The map $t\mapsto \ln(1+t)$ is continuous
at point $0$. Thus, since $\big|\frac{y}{x}\big|<1$, we have
$\lim_{p\longrightarrow
+\infty}\ln\left(\big(1+\big(\frac{y}{x}\big)^{2p+1}\big)^{\frac{1}{2p+1}}\right)=0.$
Consequently,
$$\lim_{p\longrightarrow
+\infty}\big(x^{2p+1}+y^{2p+1}\big)^{\frac{1}{2p+1}}=x\,\Big[\lim_{p\longrightarrow
+\infty}\Big(1+\big(\frac{y}{x}\big)^{2p+1}\big)^{\frac{1}{2p+1}}\Big]=x.\quad
\Box$$

Let $(M,\boxplus)$ be a magma or groupoid  that is a set $M$
equipped with a single closed binary operation $M \times M
\rightarrow M$ defined by $(x,y)\rightarrow x\boxplus y$. $M$ is {
unital }if it has a neutral element $0$. This binary operation is {
idempotent} if for all $x\in M$ $x\boxplus x=x$. It is { associative
}if for all $(x,y,z)\in M^3$, one has $(x\boxplus y)\boxplus
z=x\boxplus (y\boxplus z)=x\boxplus y\boxplus z$. $M$ has { inverse
elements} if for all $x\in M$ there is some $-x$ in $M$ such that
$x\boxplus ( -x)=(- x)\boxplus x=0$. We say that $x$ and $-x$ are {
symmetrical} if $- x$ is the inverse of $x$ and conversely. It is a
standard fact that a nontrivial group is not idempotent. In general
idempotence is compatible with  a semigroup structure that is an
algebraic structure consisting of a set together with an associative
binary operation. A semigroup generalizes  a group  to a type where
every element did not have to have a symmetrical element. In the
following a non-associative and idempotent magma is considered. The
real set $\mathbb R$ is endowed with the binary operation $\boxplus
:\mathbb R\times \Real\rightarrow \Real$ defined by

\begin{equation}\label{base}x\boxplus y=\lim_{p\longrightarrow
+\infty}x\stackrel{p}{+}y= \lim_{p\longrightarrow
+\infty}\big(x^{2p+1}+y^{2p+1}\big)^{\frac{1}{2p+1}}\end{equation}

The map $ (x,y)\mapsto x\boxplus y$ is not continuous over
$\mathbb R^2$. Moreover, this operation is not associative. For
example, one has $(1\boxplus 1)\boxplus (-1)=1\boxplus (-1)=0$ and
$1\boxplus \big( 1\boxplus (-1)\big)=1\boxplus 0=1$ which
contradicts associativity.  Associativity is replaced with a
weakened assumption which only requires that associativity works
for any pair of non symmetrical elements. By definition, for all
$x,y\in \Real_+$ one has $x\boxplus y=\max\{x,y\}$. Moreover, if
$x,y\in \Real_-$ then $x\boxplus y=\min\{x,y\}$. It follows that
the operation $\boxplus$ defines a total order on $\mathbb R_+$
and $\mathbb R_-$, but not on the whole real set. In the remainder
$(\mathbb R,\boxplus)$ denotes the set $\mathbb R$ equipped with
the operation $\boxplus$. The proof of the following lemma is left
to the reader.

\begin{prop} \label{IdG} The set $\Real$ equipped with the  operation
$\boxplus$  and the scalar multiplication satisfies the following
properties:

\noindent $(a)$ For all $x\in \Real$, $x\boxplus x=x$
$(\text{idempotence})$.

\noindent $(b)$  For all $x\in \Real$,  $x\boxplus 0=0\boxplus x=x$
$(\text{neutral element})$.

\noindent $(c)$ For all $x\in \Real$, there exists a uniqueness
symmetrical element $-x\in \Real$ such that $x\boxplus (- x)=(-
x)\boxplus x=0$ $(\text{symmetrical element})$.

\noindent $(d)$ For all $x,y\in \Real$, we have $x\boxplus
y=y\boxplus x$ $(\text{commutativity})$.

\noindent $(e)$ For all $(x,y,z)\in \Real^3$, if $x$, $y$ and $z$
are mutually non symmetrical then    $(x\boxplus y)\boxplus
z=x\boxplus (y\boxplus z)$ $(\text{weakened form of
associativity})$.

\noindent $(f)$ For all $(x,y,z)\in \Real^3$, one has   $z
(x\boxplus y)=(x\boxplus y) z=(z x)\boxplus (z y)$ $(
\text{distributivity} )$.

\end{prop}

The properties above show that $(\mathbb R,\boxplus,\cdot)$ is
endowed with some kind of ``scalar field like" algebraic
structure. It is not a scalar field because $(\mathbb R,\boxplus)$
is not a group. The next statement is an immediate consequence of
Lemma \ref{two}.

\begin{lem}\label{compar}For all $x,y\in \Real$, the following inequalities are equivalent: $(a)$
$x\leq y $; $(b)$ $ 0\leq (-x)\boxplus y $; $(c)$ $x\boxplus
(-y)\leq 0.$

\end{lem}
{\bf Proof:} First, note that the distributivity of the scalar
multiplication on the operation $\boxplus$ implies that $(b)$ and
$(c)$ are equivalent. All we need to prove is that $(a)$ and $(b)$
are equivalent. Let us prove the first implication. If $x\leq y$,
then for all natural number $p\geq 1$, we have
$((-x)^{2p+1}+y^{2p+1})^{\frac{1}{2p+1}}\geq 0$. It follows that
$\lim_{p\longrightarrow
+\infty}((-x)^{2p+1}+y^{2p+1})^{\frac{1}{2p+1}}= (-x)\boxplus
y\geq 0$, which proves $(b)$. Conversely, suppose that $(b)$
holds. By hypothesis $ 0\leq (-x)\boxplus y $. If $(-x)\boxplus
y=0$ then one has $x=y$ and $(a)$ is immediate. Suppose now that
$(-x)\boxplus y=y \geq 0$. This implies that one has $|y|\geq|-x|
= |x|$. Since $y\geq 0$, we have $y\geq x$, and we deduce
condition $(a)$. Finally, if
 $(-x)\boxplus y=-x\geq 0$, from the distributivity of the scalar
multiplication on the operation $\boxplus$, we have $x\boxplus
(-y)=x\leq 0$. Since $|x|\geq |y|$, this implies that
$x\leq y$, which ends the proof. $\Box$\\

\subsection{ Construction of a $n$-ary Operation}

In the following it is established that, though the operation
$\boxplus$ does not satisfy associativity,  it can be extended by
constructing a non-associative algebraic structure which returns to
a given $n$-tuple a real value. For all $x\in \mathbb R^n$ and all
subset $I$ of $[n]$, let us  consider the map $\xi_I[x]:\Real
\longrightarrow \mathbb Z$ defined for all $\alpha \in \Real$ by
\begin{equation}\label{defxi}\xi_I[x](\alpha)= \Card
\{i\in I: x_i=\alpha\}- \Card \{i\in I: x_i=-\alpha\}.\end{equation}
This map measures the symmetry of the occurrences of a given value
$\alpha $ in the components of a vector $x$. This map satisfies the
following properties whose the proofs are obvious and left to the
reader.

\begin{lem}\label{propxi} For all $x\in \Real^n$ and for all subset $I$ of $[n]$ the map $\xi_I[x]$  defined in
\eqref{defxi} satisfies the following properties:\\
$(a)$ $\xi_I[x]$ is an impair map, that is for all $\alpha\in
\Real$,
$\xi_I[x](-\alpha)=-\xi_I[x](\alpha)$.\\
$(b)$ For all $\alpha\in \Real$ the
map $x\mapsto \xi_I[x](\alpha)$ is impair.\\
$(c)$ If $\{\alpha,-\alpha\}\cap\{x_i:i\in I\}=\emptyset$ then
$\xi_{I}(\alpha)=0$.\\ $(d)$ For all $i\in I$,
$\xi_{I\backslash\{i\}}(x_i)=\xi_I(x_i)-1$.\\
$(e)$  If $\xi_I[x]\big(\max_{i\in I}|x_i|\big)>0$ then $\max_{i\in
I} x_i=\max_{i\in I}|x_i|$.\\$(f)$ If $\xi_I[x]\big(\max_{i\in
I}|x_i|\big)<0$ then $\min_{i\in I} x_i=-\max_{i\in I}|x_i|$.\\
$(g)$ For all subsets $I$ and $J$ of $[n]$ and all $\alpha\in
\Real$, $\xi_{I\cup
J}[x](\alpha)=\xi_{I}[x](\alpha)+\xi_{J}[x](\alpha)-\xi_{I\cap
J}[x](\alpha)$.\\
$(h)$ $\xi_\emptyset [x](\alpha)=0$, for all $\alpha \in \Real$.
\end{lem}

For all $x\in \mathbb R^n$ let $\mathcal J_I(x)$ be a subset of $I$
defined by
\begin{equation}\mathcal J_I(x)=\Big\{j\in I: \xi_I[x](x_j)\not=0\Big\}=I\backslash \big (\xi_I[x]^{-1}(0)\big).\end{equation}
$\mathcal J_{I} (x)$ is called {\bf the residual  index set } of
$x$. It is obtained by dropping from $I$ all the $i$'s such that
$\Card \{j\in I: x_j=x_i\}= \Card \{j\in I: x_j=-x_i\}$.

\begin{defn}\label{nary}For all positive natural number $n$ and
for all subset $I$ of $[n]$, let $\digamma_I: \Real^n
\longrightarrow \Real $ be the map defined for all $x\in \Real^n$ by

\begin{equation}\digamma_{ I}(x)=\left\{\begin{matrix}\max_{i\in \mathcal
J_I (x)}x_{i} &\text{ if }&\xi_I[x]\big (\max_{i\in
\mathcal J_I(x)}|x_i|\big)>0\\
\min_{i\in \mathcal J_I(x)}x_i &\hbox{ if }&\xi_I[x]\big(\max_{i\in \mathcal J_I(x)}|x_i|\big)<0\\
 0 &\text{ if }&\xi_I[x]\big (\max_{i\in
\mathcal J_I(x)}|x_i|\big)=0.
\end{matrix}\right.\end{equation}
where $\xi_I[x]$ is the map defined in \eqref{defxi} and $\mathcal
J_I(x)$ is the residual index set of $x$. The operation that takes
an $n$-tuple $(x_1,....,x_n)$ of $\Real^n$ and returns a single real
element $\digamma_I(x_1,...,x_n)$ is called a $n$-ary extension of
the binary operation $\boxplus.$
\end{defn}

Notice that, if $\mathcal J_I(x)=\emptyset$ if and only if
$\xi_I[x](\,\max_{i\in \mathcal J_I(x)}|x_i|\,)=0$.   To see the key
idea of the definition above let us define the generalized sum of
$n$ real numbers $x_1,...,x_n$ as $ S_p(x_1,...,x_n)=\Big(\sum_{i\in
[n]}{x_i}^p\Big)^{\frac{1}{p}}$, say a H\"older sum. When one
consider the subsequence of  pair natural numbers, this generalized
sum has the limit: $ \lim_{p\longrightarrow
+\infty}S_{2p}(x_1,...,x_n)=\lim_{p\longrightarrow
+\infty}\Big(\sum_{i\in
[n]}{x_i}^{2p}\Big)^{\frac{1}{2p}}=\max_{i\in [n]}|x_i|. $
 The case where the generalized sum is defined with respect to
 the
  impair natural numbers is analyzed in this section. It is shown below  that $\digamma_{I}(x)$ is the limit of the
generalized sum $S_{2p+1}(x_1,...,x_n)$.

\begin{prop} \label{LimMaxFunc}
For all natural number $n\geq 1$ and all $x\in \Real^{n}$, if $I$ is
a nonempty subset of $[n]$ then:
\begin{align*}\digamma_{I}(x)=\lim_{p\longrightarrow \infty}\Big( \sum_{i\in
I}{x_i}^{2p+1}\Big)^{\frac{1}{2p+1}}=\lim_{p\longrightarrow
\infty}\stackrel{\varphi_p}{\sum_{i\in I}}x_i.
\end{align*}

\end{prop}

\noindent {\bf Proof:} Let $\mathcal J_{I}(x)$ be the residual index
set of $x$. We have

$$\Big(\sum_{i\in I}{x_{i}}^{2p+1}\Big)^{\frac{1}{2p+1}}=\Big(\sum_{i\in
  I\backslash \mathcal J_{I}(x)}{x_{i}}^{2p+1}+\sum_{i\in
\mathcal J_{I}(x)}{x_{i}}^{2p+1}\Big)^{\frac{1}{2p+1}}.
$$ By definition, there exists a partition of $ I\backslash \mathcal J_{I}(x)$ whose any block contains two symmetric elements. Hence it follows
that
$$\sum_{i\in
I\backslash \mathcal J_{I}(x)}{x_{i}}^{2p+1}=-\sum_{i\in I\backslash
\mathcal J_{I}(x)}{x_{i}}^{2p+1}=0.$$ Hence, we deduce that $$
\Big(\sum_{i\in I}{x_{i}}^{2p+1}\Big)^{\frac{1}{2p+1}}
=\Big(\sum_{i\in \mathcal
J_{I}(x)}{x_{i}}^{2p+1}\Big)^{\frac{1}{2p+1}}.$$ Suppose that
$\mathcal J_{I}(x)=\emptyset$. In such a case
$$\Big(\sum_{i\in I}{x_{i}}^{2p+1}\Big)^{\frac{1}{2p+1}}
=\Big(\sum_{i\in \mathcal
J_{I}(x)}{x_{i}}^{2p+1}\Big)^{\frac{1}{2p+1}}=0.$$ Consequently
$\lim_{p\longrightarrow \infty}\Big(\sum_{i\in
I}{x_{i}}^{2p+1}\Big)^{\frac{1}{2p+1}}=0=\digamma_I(x)$ which proves
this case.

Suppose now that $\mathcal J_{I}(x)\not=\emptyset$. Let us denote
$$\mathcal M_I(x)=\{i\in \mathcal J_I(x): |x_i|=\max_{i\in \mathcal
J_I(x)}|x_i|\}.$$ Then, from the definition of map $\xi_I[x]$ in
equation \eqref{defxi}, we have
$$\sum_{i\in
\mathcal M_I(x)}\Big(\frac{x_{i}}{\max_{i\in \mathcal
J_I(x)}|x_i|}\Big)^{2p+1}=\xi_I[x]\big(\max_{i\in \mathcal
J_I(x)}|x_i|\big) .$$

It follows that for all $x\in \Real^n$,  {\small
\begin{align} \Big(\sum_{i\in
\mathcal J_{I}(x)}{x_{i}}^{2p+1}\Big)^{\frac{1}{2p+1}} &=\max_{i\in
\mathcal J_{I}(x)} |x_{i}|\Big(\sum_{i\in \mathcal J_{I}(x)}
\frac{{x_{i}}^{2p+1}}{\max_{i\in \mathcal J_{I}(x)}|x_{i}|^{2p+1}}\Big)^{\frac{1}{2p+1}}\label{factor1}\\
&= \left(\max_{i\in \mathcal J_{I}(x)}|x_{i}|\right)\alpha_p
(x).\nonumber
\end{align}}
where
$$\alpha_p (x)=\Big(\xi_I[x]\big (\max_{i\in
\mathcal J_I(x)}|x_i|\big ) + \sum_{i\notin \mathcal M_{ I}(x)}
\big(\frac{x_{i}}{\max_{i\in \mathcal J_{I}(x)}
|x_{i}|}\big)^{2p+1}\Big)^{\frac{1}{2p+1}}.$$

\noindent We need to compute the limit of $\alpha_p (x)$, when
$p\longrightarrow \infty$. Clearly, for all $i\notin \mathcal
M_{I}^{}(x)$, we have
\begin{equation}\label{major1}
\Big|\frac{x_{i}}{\max_{i\in \mathcal
J_{I}(x)}|x_{i}|}\Big|<1.\end{equation} Since $\mathcal
J_I(x)\not=\emptyset$, $\xi_I[x]\left(\max_{i\in
\mathcal J_I(x)}|x_i|\,\right)\not=0$, hence we consider two cases:\\

\noindent $(i)$ $\xi_I[x]\left (\max_{i\in
\mathcal J_I(x)}|x_i|\right )>0$. \\

For the sake of simplicity, define $a=\xi_I[x](\max_{i\in \mathcal
J_I(x)}|x_i|\,)$ and $b_i= |\frac{x_{i}}{\max_{i\in \mathcal
J_{I}(x)}|x_{i}|}|$ for each $i$. We then obtain $\alpha_p(x)=
\big(a + \sum_{i\notin \mathcal M_{ {I}}(|x|)}
b_i^{2p+1}\big)^{\frac{1}{2p+1}}.$ Moreover, from \eqref{major1}, we
have $|b_i|<1$ for all $i\notin \mathcal M_{ I}^{}(|x|)$. By
hypothesis $a>0$ and we deduce that $\lim_{p\longrightarrow +\infty}
\ln (a + \sum_{i\notin \mathcal M_{ {I}}(|x|)} b_i^{2p+1})=\ln (a)$.
Hence, we have
$$\lim_{p\longrightarrow +\infty}\ln \alpha_p(x)=
\lim_{p\longrightarrow +\infty}\frac{\ln (a + \sum_{i\notin \mathcal
M_{ {I}}(|x|)} b_i^{2p+1})}{2p+1}=0.$$ Thus, $\lim_{p\longrightarrow
+\infty} \alpha_p(x)=1$. Hence, from \eqref{factor1}, we deduce
that:
$$\lim_{p\longrightarrow \infty}\Big(
\sum_{i\in \mathcal
J_{I}(x)}{x_i}^{2p+1}\Big)^{\frac{1}{2p+1}}=\max_{i\in \mathcal
J_{I}(x)}|x_i|=\digamma_I(x).$$

\noindent $(ii)$ $\xi_I[x]\left (\max_{i\in
\mathcal J_I(x)}|x_i|\right)<0$. \\

 Applying $(i)$, we then
obtain
$$-\lim_{p\longrightarrow \infty}\Big( \sum_{i\in
\mathcal
J_{I}(x)}{x_i}^{2p+1}\Big)^{\frac{1}{2p+1}}=\lim_{p\longrightarrow
\infty}\Big( \sum_{i\in \mathcal
J_{I}(x)}{(-x_i)}^{2p+1}\Big)^{\frac{1}{2p+1}}=\max_{i\in \mathcal
J_{I}(x)}|x_i|.$$ Thus  $$\lim_{p\longrightarrow \infty}\Big(
\sum_{i\in \mathcal
J_{I}(x)}{x_i}^{2p+1}\Big)^{\frac{1}{2p+1}}=-\max_{i\in \mathcal
J_{I}(x)}|x_i|=\min_{i\in \mathcal J_{I}(x)}x_i=\digamma_I(x).\quad
\Box$$

Let us introduce for all $n$-tuple $x=(x_1,...,x_n)$ the operation
defined by:

\begin{equation}\label{defrecOp}\bigboxplus_{i\in
I}x_i=\lim_{p\longrightarrow \infty}\Big( \sum_{i\in
I}{x_i}^{2p+1}\Big)^{\frac{1}{2p+1}}=\lim_{p\longrightarrow
\infty}\stackrel{\varphi_p}{\sum_{i\in I}}x_i.
\end{equation}

Clearly, this operation encompasses as a special case the binary
operation defined in equation \eqref{base}. From Fact \ref{two} and
Definition \ref{nary} if $n=2$ and $I=\{1,2\}$,  then, for all
$(x_1,x_2)\in \Real^2$
\begin{equation*}\bigboxplus_{i\in \{1,2\}}x_i=x_1\boxplus
x_2.\end{equation*}

\begin{expl}
Suppose that $x=(2,3,-2,-3,\frac{3}{2},-3,3,-\frac{1}{2})$. First,
note that $[8]=\{1,\cdots ,8\}$ and $\{x_i: i\in
[8]\}=\{1,-2,2,-3,3,-\frac{1}{2}\}$. We have $\{i: x_i=3\}=\{2,7\}$
and $\{i: x_i=-3\}=\{4,6\}$. Therefore, $\Card\{i: x_i=3\}=\Card
\{i: x_i=-3\}=2$, $\xi_{[8]}[x](x_i)=0$ for $i\in \{2,4,6,7\}$.
Moreover: $\{i: x_i=2\}=\{1\}$ and $\{i: x_i=-2\}=\{3\}$.
Consequently, $\Card \{i: x_i=2\}=\Card \{i: x_i=-2\}=1$ and
$\xi_{[8]}[x](x_i)=0$ for $i\in \{1,3\}$. Hence, we have
$$\mathcal J_{[8]}(x)=[8]\backslash \big(\{2,4,6,7\}\cup\{1,3\}\big)=\{5,8\}.$$
Therefore $\mathcal J_{[8]}(x)=\{5,8\}$. Hence $$\bigboxplus_{i\in
[8]}x_i=\bigboxplus_{i\in \{5,8\}}x_i=\frac{3}{2}\boxplus
(-\frac{1}{2})=x_5=\frac{3}{2}.$$
\end{expl}

\subsection{Some Algebraic Properties} A few immediate properties whose the
proofs are obvious are established in the next Lemma.

\begin{prop}  \label{propalg}For all $x\in \Real^n$ and all nonempty subset $I$ of $[n]$, we
have:\\
$(a)$ If $\mathcal J_I(x)\not=\emptyset$ then there is some $i_0\in
I$ such that $x_{i_0}=\bigboxplus_{i\in I}x_i$. Moreover, $\xi[x](x_{i_0})>0$.\\
$(b)$ If $\alpha\in \Real$ and $x_i=\alpha$ for all $i\in I$ then
$\bigboxplus_{i\in I}x_i=\alpha$.\\ $(c)$ Moreover, if all the
elements of the family $\{x_i\}_{i\in I}$ are mutually non
symmetrical, then: $\bigboxplus_{i\in
I}x_i=\arg\max_{x_i}\{|x_i|,i\in I\}.$\\ $(d)$ For all $\alpha\in
\Real$, one has:
\begin{equation*}
\alpha \Big(\bigboxplus_{i\in I}x_i\Big)= \bigboxplus_{i\in
I}(\alpha x_i).
\end{equation*}
$(e)$ Suppose that $x\in\epsilon \Real_+^n$ where $\epsilon$ is
$+1$ or $-1$. Then $\bigboxplus_{i\in I}x_i=\epsilon \max_{i\in I}
\{\epsilon x_i\}$.\\
$(f)$ We have $|\bigboxplus_{i\in I} x_i|\leq \bigboxplus_{i\in
I}|x_i|$.\\
$(g)$ For  all permutation $\sigma:I\rightarrow I$, we have
$\bigboxplus_{i\in I} x_i=\bigboxplus_{i\in I}x_{\sigma(i)}$.\\$(h)$
If there exists $j,k\in I$ with $j\not=k$ and $x_j+x_k=0$, then
$\bigboxplus_{i\in I\backslash\{j,k\}}x_i=\bigboxplus_{i\in I}x_i. $
Moreover $\bigboxplus_{i\in I}x_i=\bigboxplus_{i\in \mathcal
J_I(x)}x_i. $
\end{prop}
{\bf Proof:}   $(a)$  By hypothesis, the subset $J=\left\{i\in I:
\xi_{I}[x](x_i)\not=0\right\}$ is nonempty. Therefore, there exists
some $i_0\in \mathcal J_I(x)$ such that $|x_{i_0}|\geq |x_i|$ for
all $i\in \mathcal J_I(x)$. There are two possibilities. If
$\xi_I[x](\max_{i\in \mathcal J_I(x)}|x_i|)> 0$, then, from Lemma
\ref{propxi}.e, $\max_{i\in \mathcal J_I(x)}|x_i|=\max_{i\in
\mathcal J_I(x)}x_i=x_{i_0}$. Thus
 $\xi_I[x](x_{i_0})> 0$. If $\xi_I[x](\max_{i\in \mathcal J_I(x)}|x_i|)< 0$, from
 Lemma
\ref{propxi}.f, we have $\min_{i\in \mathcal J_I(x)}x_i=-\max_{i\in
\mathcal J_I(x)}|x_i|=x_{i_0}$. Since $\xi_I[x]$ is an impair map,
this also implies that $\xi_I[x](x_{i_0})> 0$. $(b)$ is immediate
setting $x_i=\alpha$ for all $i\in I$. $(c)$ If all the $x_i$'s are
mutually non symmetrical, then $\xi_I[x](x_i)\not=0$ for all $i\in
I$. Hence, $\xi_I[x](\max_{i\in I}|x_i|\,)\not=0$ and there is some
$\epsilon\in \{-1,1\}$ such that $\digamma_I(x)=\epsilon \max_{j\in
\mathcal J_I(x)} |x_i|$, we deduce $(c)$. $(d)$ Since the scalar
multiplication is distributive on addition, it is an immediate
consequence of Proposition \ref{LimMaxFunc}. $(e)$ If $x\in
\Real_+^n$ then $\bigboxplus_{i\in I}x_i=\max_{i\in I}x_i$. $x\in
-\Real_+^n$ implies that $\xi_I[x](\max_{i\in I}|x_i|\,)<0$ and
$\bigboxplus_{i\in I}x_i=\min_{i\in I}x_i=-\max_{i\in I}(-x_i)$.
$(f)$ For all $x\in \mathbb R^n$ and all $I\subset [n]$, we have
$\digamma_I (|x|)=\max_{i\in I}|x_i|$. Moreover, by definition there
is some $\delta\in \{-1,0,1\}$ such that $\digamma_I(x)=\delta
\max_{j\in \mathcal J_I(x)}\delta x_i$. Therefore $\digamma_I
(|x|)\geq |\digamma_I (x)|$. $(g)$ Given a nonempty subset $I$ of
$[n]$, the H\"older sum is independent of any permutation of the
index set $I$. Therefore, from Proposition \ref{LimMaxFunc},
 we deduce $(g)$. $(h)$ In such a case
 $\{i\in I\backslash\{j,k\}:\xi_I[x](x_i)\not=0\}=
  \{i\in I: \xi_I[x](x_i)\not=0\}\}$. Therefore $\mathcal J_{I\backslash\{j,k\}}(x)=\mathcal J_{I}(x)$, which proves the first part of the
  statement. Since $\mathcal J_{\mathcal J_I(x)}(x)=\mathcal J_I(x)$, the second part is immediate. $\Box$. \\

 In the following we introduce the operation $ \langle \cdot,\cdot\rangle_\infty :\Real^n\times \Real^n\longrightarrow
 \Real$ defined for all $x,y\in \Real^n$ by $\langle x,y\rangle_\infty =\bigboxplus_{i\in
 [n]}x_iy_i$. Let $\|\cdot\|_\infty$ be the Tchebychev
 norm defined by $\|x\|_\infty=\max_{i\in [n]}|x_i|$. The next
 result is an immediate consequence of  Proposition \ref{propalg}.

 \begin{prop}For all $x,y\in \Real^n$, we have:\\
 $(a)$  $\sqrt{\langle x,x\rangle_\infty} =\|x\|_\infty$.\\
  $(b)$  $|\langle x,y\rangle_\infty| \leq \|x\|_\infty \|y\|_\infty$.\\
   $(c)$ For all $\alpha\in \Real$, $\alpha \langle x,y\rangle_\infty= \langle \alpha x,y\rangle_\infty=\langle  x,\alpha y\rangle_\infty$.
\end{prop}
{\bf Proof:} $(a)$ By definition $\langle
x,x\rangle_\infty=\bigboxplus_{i\in [n]}{x_i}^2=\max_{i\in
[n]}{x_i}^2={\|x\|_\infty}^2$, which ends the proof. $(b)$ From
Proposition \ref{propalg}.f, $|\langle
x,y\rangle_\infty|=|\bigboxplus_{i\in [n]}x_iy_i|\leq
\bigboxplus_{i\in [n]}|x_iy_i|\leq (\max_{i\in [n]}|x_i|)(\max_{i\in
[n]}|y_i|)$, which proves $(b)$. $(c)$ is immediate from Proposition
\ref{propalg}.d. $\Box$\\

The next statement establishes a key  property resulting from
Proposition \ref{preassoc}.

\begin{prop} \label{preassoc}Suppose that $x=(x_1,...,x_n)\in \Real^n$. For all
nonempty subset $I$ of $[n]$ and all $i\in I$:
$$\Big[x_i \boxplus
\big(\bigboxplus_{j\in I\backslash \{i\}}x_j\big)\Big]\in \Big
\{0, \bigboxplus_{j\in I}x_j\}$$ and
$$\bigboxplus_{i\in
I}x_i=\bigboxplus_{i\in I}\Big[x_i \boxplus \big(\bigboxplus_{j\in
I\backslash \{i\}}x_j\big)\Big].$$

\end{prop}
{\bf Proof:} If $\mathcal J_I(x)=\emptyset$, then this property is
immediate. In such case, since $I$ is nonempty, from Proposition
\ref{propalg}.a, there is some $i_0\in I\backslash\{i\}$ such that
$x_{i_0}=-x_i$. Therefore $x_i \boxplus \big(\bigboxplus_{j\in
I\backslash \{i\}}x_j\big)=0$. Suppose now that  $\mathcal
J_{I\backslash \{i\}}(x)\not=\emptyset$ and let us consider four
cases:

 $(i)$  $|x_i|>|\bigboxplus_{j\in I}x_j|$. This implies that $i\in I\backslash \mathcal J_I(x) $. Thus $\xi_I[x](x_i)=0$. Hence,
 $\xi_{I\backslash
\{{i}\}}[x](x_i)<0$ and from Proposition \ref{propalg}.a
$\bigboxplus_{j\in I\backslash \{{i}\}}x_j\not=x_i$. Moreover
$\xi_{I\backslash \{{i}\}}[x](-x_i)>0$, and by hypothesis,
$|x_i|\geq |x_j|$ for all $j\in \mathcal J_{I\backslash \{i\}}(x)$.
Therefore $\bigboxplus_{j\in I\backslash \{{i}\}}x_j=-x_i.$ It
follows that
$$x_{{i}} \boxplus \big(\bigboxplus_{j\in I\backslash
\{{i}\}}x_j\big)=x_{{i}}\boxplus(-x_{{i}})=0,$$ which proves this
case.

$(ii)$  $x_i=\bigboxplus_{j\in I}x_j$. By definition, this implies
that
  $\xi_I[x](x_i)>0$. Since
  $\xi_{I\backslash\{i\}}[x](x_i)=\xi_{I}[x](x_i)-1$, one has $\xi_{I\backslash \{i\}}[x](x_i)\geq 0$ and,
consequently, $\xi_{I\backslash \{i\}}[x](-x_i)\leq 0$. Thus, from
Proposition \ref{propalg}.a, $\bigboxplus_{j\in I\backslash
\{i\}}x_j \not=-x_i$. Moreover, $|\bigboxplus_{j\in I\backslash
\{i\}}x_j |\leq |x_i|$ and we have
$$x_{{i}} \boxplus \big(\bigboxplus_{j\in I\backslash
\{{i}\}}x_j\big)=x_{{i}}=\bigboxplus_{j\in I}x_j.$$

$(iii)$  $x_i=-\bigboxplus_{j\in I}x_j$. Equivalently, we have
$-x_i=\bigboxplus_{j\in I}x_j$ and, from Proposition
\ref{propalg}.a, this implies that there is some $i_0\in I$ such
that $x_{i_0}=\bigboxplus_{j\in I}x_j=-x_i$ with
 $\xi_I[x](x_{i_0})>0$.  Since $\mathcal J_I(x)\not=\emptyset$, $x_{i_0}\not=0$ and $x_{i_0}\not=x_i$. It follows
 that ${i_0}\not=i$. Thus
  $\xi_{I\backslash \{i\}}[x](x_{i_0})>0$. Therefore, $\bigboxplus_{j\in I\backslash
 \{i\}}x_j=x_{i_0}=-x_i$.
Consequently:
$$x_{{i}} \boxplus \big(\bigboxplus_{j\in I\backslash
\{{i}\}}x_j\big)=x_{{i}} \boxplus (-x_i)=0.$$

 $(iv)$  $|x_i|<|\bigboxplus_{j\in I}x_j|$.
 Moreover, from Proposition \ref{propalg}.a there is some $i_0\in I$ such that $x_{i_0}=\bigboxplus_{j\in I}x_j$. This implies that
 $ \xi_I[x](x_{i_0})>0$. Therefore,
since $|x_{i_0}|>|x_i|$, we have $\xi_{I\backslash
\{i\}}[x](x_{i_0})>0$ and it follows that
$|x_i|<|\bigboxplus_{j\in I\backslash
 \{i\}}x_i|$. Hence $x_{{i}} \boxplus \big(\bigboxplus_{j\in I\backslash
\{{i}\}}x_j\big)=x_{{i_0}}=\bigboxplus_{j\in I}x_j,$ which ends the
proof of the first part of the
statement.\\

 To prove the second part of the statement, we need to establish that there exists some
$i\in I$ such that $x_{i}\boxplus \big(\boxplus_{i\in
I\backslash\{i\}}x_i\big)=x_i$. If $\mathcal J_I(x)=\emptyset$ then
$\bigboxplus_{j\in I}x_j=0$ and from the statement above $\Big[x_i
\boxplus \big(\bigboxplus_{j\in I\backslash \{i\}}x_j\big)\Big]\in
\{0\}$ for all $i$. In such a case, this property is obviously true.
Suppose that $\mathcal J_I(x)\not=\emptyset$. Recall that from
Proposition \ref{propalg}.a there is some $i\in I$ such that
$x_{i}=\boxplus_{j\in I}x_j$. Then, using
$(ii)$, the second statement follows. $\Box$\\

For example, for all  $x,y,z\in \Real.$ we have the identities: $
x\boxplus y=x\boxplus y\boxplus  x\boxplus y$ and $ x\boxplus
y\boxplus z=\big [x\boxplus (y\boxplus  z)\big]\boxplus
\big[y\boxplus (z\boxplus  x)\big ] \boxplus \big [z\boxplus
(x\boxplus y)\big ].$

 Let  $\Lambda:
\Real^n\longrightarrow \Real^n$ be the map defined by:
\begin{equation}
\Lambda(x_1,...,x_n)= \Big(x_1 \boxplus \big(\bigboxplus_{j\in
[n]\backslash \{1\}}x_j\big),...,x_n \boxplus
\big(\bigboxplus_{j\in [n]\backslash \{n\}}x_j\big)\Big).
\end{equation}

\begin{lem}
\label{assoc} Suppose that $x=(x_1,...,x_n)\in \Real^n$. Let
$\epsilon\in \{-1,1\}$ such that $\bigboxplus_{i\in I}x_i\in
\epsilon \Real_+$. Then:\\
$(a)$ $\bigboxplus_{i\in I}x_i =\epsilon \max _{i\in I}\{\epsilon
\Lambda_i(x)\};$ \\ $(b)$ For all $i,j,k\in [n]$, we have:
$\Lambda_i(x)\boxplus \Lambda_j(x) \boxplus
\Lambda_k(x)=\Big(\Lambda_i(x)\boxplus \Lambda_j(x) \Big)\boxplus
\Lambda_k(x)=\Lambda_i(x)\boxplus \Big(\Lambda_j(x) \boxplus
\Lambda_k(x)\Big).$
\end{lem}
{\bf Proof:} $(a)$ From Proposition \ref{preassoc} we have
$\Lambda_i(x)\in  \{0,\bigboxplus_{i\in I}x_i \}$ for all $i\in
I$. Hence  $\Lambda_i(x)\in \epsilon \Real_+$ for all $i$ and from
Proposition \ref{preassoc}.c, the result follows. $(b)$ is an
immediate consequence of
$(a)$.  $\Box$\\

\begin{expl}
Let $x=(4, -3,-4, 2,3, 2,-2)\in \Real^7$. We have $\mathcal
J_{[7]}(x)=\{4,6,7\}$ and $[7]\backslash \mathcal
J_{[7]}(x)=\{1,2,3,5\}$. We have $4\boxplus (-3)\boxplus
(-4)\boxplus 2\boxplus 3 \boxplus 2 \boxplus (-2) =2$. Moreover,
$\Lambda_1(x)=4\boxplus \big( (-3)\boxplus (-4)\boxplus 2\boxplus
3 \boxplus 2 \boxplus (-2)  \big)=4 \boxplus (-4)=0$; Similarly we
obtain $\Lambda_2(x)=(-3) \boxplus 3=0$; $\Lambda_2(x)=(-4)
\boxplus 4=0$; $\Lambda_4(x)=2 \boxplus 0=2$; $\Lambda_5(x)=3
\boxplus (-3)=0$; $\Lambda_6(x)=2 \boxplus 0=2$;
$\Lambda_7(x)=(-2) \boxplus 2=0$. It follows that
$$\Lambda(x)=(0,0,0,2,0,2,0)=\big(0,0,0,\bigboxplus_{i\in [7]}x_i,0,\bigboxplus_{i\in [7]}x_i,0\big).$$
\end{expl}

\begin{lem}\label{symcopos}
Let $n$ be a positive natural number and $I$ be a nonempty subset of
$[n]$. Let $\mathfrak P(I)=\{I_j:j\in [m]\}$ be a partition of $I$
with $m$ nonempty subsets $I_j$.  If for all $(j,k)\in [m]\times
[m]$
$$\Big(\bigboxplus_{i\in I_j} x_{i}\Big)+\Big(\bigboxplus_{i\in I_k} x_{i}\Big)\not=0\;\text{ then }
\;\bigboxplus_{j\in [m]}\Big(\bigboxplus_{i\in I_j}
x_{i}\Big)=\bigboxplus_{i\in I}x_{i}.$$

\end{lem}
{\bf Proof:} For all $j\in [m]$, let us denote
$y_j=\bigboxplus_{i\in I_j} x_{i}$.  By hypothesis  the $y_j$'s are
mutually non symmetrical, it follows that there exists some $j_0\in
[m]$ such that $$y_{j_0}=\bigboxplus_{j\in [m]}
y_j=\arg\max_{y_j}\{|y_j|:j\in [m]\}.$$ Therefore, for all $i\in I$
such that $|x_i|>|y_{j_0}|$, and all $j\in [m]$
$\xi_{I_j}\,[x](x_i)=0$. However, by hypothesis $\mathfrak P(I)$ is
a partition of $I$. Hence $I=\bigcup_{j\in [m]}I_j$  with $I_j\cap
I_k=\emptyset$, for all $j\not=k$. Thus, for all $i\in I$ such that
$|x_i|>|y_{j_0}|$ we have from Lemma \ref{propxi}.g
$\xi_{I}[x](x_i)=\sum_{j\in [m]}\xi_{I_j}[x](x_i)=0$. It follows
that
$$\mathcal J_I(x)\subset \{i\in I: |x_{i}|\leq |y_{j_0}|\}.$$

Therefore, $|\bigboxplus_{j\in [m]}\Big(\bigboxplus_{i\in I_j}
x_{i}\Big)|=|\bigboxplus_{j\in [m]}y_j|=|y_{j_0}|\geq
|\bigboxplus_{i\in I}x_{i}|.$ From Proposition \ref{propalg}.a,
$y_{j_0}\not=0$ implies that there is some $i_0\in I_{j_0}$ and such
that $x_{i_0}=\bigboxplus_{i\in I_{j_0}}x_i=y_{j_0}$. Since
$|x_i|\leq |x_{i_0}|$ for all $i\in \mathcal J_I(x)$ and the $y_j$'s
are not symmetrical it follows that
$\xi_I[x](x_{i_0})>0$ which ends the proof. $\Box$\\

\subsection{Euclidean Orthant, Absolute Value and Upper Semi-Lattice
Structure}\label{orthan}

The algebraic structure $(\Real,\boxplus,\cdot)$  can be extended
to $\Real^n$. Suppose that $x,y\in \Real^n$, and let us denote
\begin{equation}
x\boxplus y=(x_1\boxplus y_1,\cdots,x_n\boxplus y_n).
\end{equation}
Moreover, let us consider $m$ vectors $x_1,...,x_m\in \Real^n$, and
define
\begin{align}
\bigboxplus_{j\in [m]}x_j&=\Big(\bigboxplus_{j\in [m]}
x_{j,1},\cdots,\bigboxplus_{j\in [m]} x_{j,n}\Big).
\end{align}

Let the triple $(\Real^n,\boxplus,\cdot)$ denotes the
$n$-dimensional Euclidean vector space  equipped with the operation
binary operation $(x,y)\mapsto x\boxplus y$ and the external scalar
multiplication of vectors by real numbers $\cdot$.

For all $(x,y)\in \Real_+^n\times \Real_+^n$ we have $x\boxplus
y=x\vee y$. Moreover, for all  $(x,y)\in\Real_-^n\times
\Real_-^n$, $x\boxplus y=x\wedge y$ where $\vee$ and $\wedge$
respectively denote the maximum and the minimum with respect to
partial order of $\mathbb R^n$ associated to the positive cone,
that is, the coordinatewise supremum and infimum. For all $x$ and
$y$ in $\mathbb R^n$, $ x\leq y$ means $y - x \in \mathbb R_+^n$.
It follows that given a subset $\{x_1,...,x_m\}$ of $\Real_+^m$,
we have \begin{equation} \bigboxplus_{j\in [m]}x_j=\bigvee_{j\in
[m]}x_j .\end{equation} If $\{x_1,...,x_m\}$ is a subset of
$\Real_-^m$ then \begin{equation} \bigboxplus_{j\in
[m]}x_j=\bigwedge_{j\in [m]}x_j .\end{equation}

For all $x,y\in \mathbb R^n$, let us denote $x\boxdot
y=(x_1y_1,...,x_ny_n)$. In  the following we say that two vectors
$x,y\in \Real^n$ are {\bf copositive} if
\begin{equation}x\boxdot y\in \Real_+^n.\end{equation}
We say that a subset $K$ of $\Real^n$ is copositive if for all
$x,y\in K$ one has $x\boxdot y\geq 0$. For all subset $L$ of
$\mathbb R^n$,  $K$ is copositive and maximal in $L$ if there does
not exists a copositive subset $K'\subset L$ which contains $K$. A
$n$-dimensional ortant in $\Real^n$ is  copositive and maximal in
$\Real^n$. Equivalently, a $n$-dimensional orthant in $\Real^n$ is a
subset defined by a system of inequalities: $\epsilon_ix_i \geq 0$
for any $i\in [n]$, where each $\epsilon_i$ is $+1$ or $-1$. A
$n$-dimensional closed orthant $K$ of $\Real^n$ can be written
$K=\prod_{i=1}^n\big(\epsilon_i\Real_+\big) $.  Let
$\Psi_K:\Real_+^n\longrightarrow K$ be the map defined by $\Psi_K
(x)=(\epsilon_1 x_1,...,\epsilon_n x_n)$ with $|\epsilon_i|=1$, for
all $i\in [n]$. $\psi_K$ is a linear homeomorphism such that
$\psi_K(\Real_+^n)=K$ and one has ${\psi_K}^{-1}=\psi_K$, which
implies that $\psi_K(K)=\Real_+^n$.

For all $x\in \Real^n$, let us denote  $|x|=(|x_1|,...,|x_n|)$. Let
$K$ be a $n$-dimensional orthant and let us consider the binary
relation defined  by $x\leqslant y\iff |x|\leq |y|$. $\geqslant$ is
a partial order over  $K$. For all $x$, $y$, and $z$ in $K$, we have
 $x \leqslant x$ (reflexivity); if $x \leqslant y$ and $y
\leqslant x$ then $x = y$ (antisymmetry); if $x \leqslant y$ and $y
\leqslant z$ then $x \leqslant z$ (transitivity). A $n$-dimensional
closed orthant $K$ equipped with the partial order $\leqslant$ is  a
partially ordered set (or a poset). Then $\boxplus$ is a join on
$K$, and the triple $(K,\boxplus,\geqslant)$ is an
upper-semilattice.

If $\{x_1,...,x_m\}$ is a subset of $K$ then
\begin{equation} \bigboxplus_{j\in [m]}x_j=\Psi_K \Big(\bigvee_{j\in
[m]}\Psi_K(x_j) \Big).\end{equation}

\section{On Some Idempotent Convex Structure }

 A subset $C$ of a $\Real_+^n$ is $\mathbb B$-convex if and only if for all
$t\in [0,1]$ and all $x,y\in C$, $x\vee ty\in C$. Equivalently, we
say that a subset $C$ of a $n$-dimensional
 orthant $K$ is $\mathbb B$-convex if and only if for all
$t\in [0,1]$ and all $x,y\in C$, $x\boxplus ty\in C$. Such a
definition is equivalent to that one proposed in \cite{bh}.  It is
also the definition proposed further in the paper to define $\mathbb
B$-convex sets on the whole Euclidean vector space. Equivalently, a
subset $C$ of $K$ is $\mathbb B$-convex if and only if $\Psi_K (C)$
is a $\mathbb B$-convex subset of $\Real_+^n$. For all finite subset
$A=\{x_1,...,x_m\}$ of $K$ the smallest $\mathbb B$-convex set which
contains it is $\mathbb B[A]=\left\{\bigboxplus_{i=1, \cdots, m}t_i
x_i : t_i\in [0,1], \max_{i \in[ m]} t_i=1\right\}$. For the sake of
simplicity, let $\mathbb B[x,y]$ denote the $\mathbb B$-convex hull
of $\{x,y\}$ for all $x,y\in K$.

The binary operation $\boxplus$ yields a simple formulation of
$\mathbb B$-convexity on each orthant.  However,  the problem to
solve is much more complex over $\Real^n$. Suppose for example that
$x,y\in \Real^n$, $|x|=|y|$ and $x\not=y$, then $\big \{t x\boxplus
s y: \max\{t,s\}=1,t,s\geq 0\}=\{x,x\boxplus y,y\big \}$ that is not
a path-connected subset of $\Real^n$.

\subsection{An  Extended Definition of $\mathbb B$-convexity}

In ~\cite{bh} $\mathbb{B}$-convexity is introduced as a limit of
linear convexities. More precisely, for all $p\in \mathbb{N}$ the
$\Phi_p$-convex hull of  a finite set $A\subset \mathbb{R}^n$ is
defined by:
\begin{equation}Co^{\Phi_p}(A)=\Big\{\stackrel{\varphi_p}{\sum_{i\in
[m]}}t_i\stackrel{\varphi_p}{.}x_i : \stackrel{\varphi_p}{\sum_{i\in
[m]}}t_i=\varphi_p^{-1}(1),
 \varphi_p(t_i)\geq 0, i\in [m]\Big\}\end{equation} which can be
rewritten:
$$Co^{\Phi_p}(A)=\left\{\Phi_p^{-1}\Big(\sum_{i\in [m]}t_i^{2p+1}{.}\Phi_p(x_i)\Big) :
\Big(\sum_{i\in [m]}t_i^{2p+1}\Big)^{\frac{1}{2p+1}}=1 , \,
t_i\geq 0, i\in [m]\right\}.$$ This is basically the approach of
Ben-Tal  \cite{ben} and Avriel  \cite{avr1}.

Equivalently, one has
$Co^{\Phi_p}(A)=\Phi_p^{-1}\Big(Co\big(\Phi_p(A)\big)\Big)$. Recall
that, for all $x,y\in \Real^n$,
$x\stackrel{p}{+}y=x\stackrel{\varphi_p}{+}y$. For simplicity,
throughout the remainder of the paper we denote for all subset $L$
of $\Real^n$ $ Co^{p}(L)=Co^{\Phi_p}(L)$.

From Briec and Horvath \cite{bh} a subset $L$ of $\Real^{n}$ is
$\mathbb{B}$-convex if for all finite subset $A\subset L$ the
$\mathbb{B}$-polytope $Co^{\infty}(A)=Ls_{p\longrightarrow \infty}
Co^{p}(L)$ is contained in $ L$. In the following, we show that the
Painlev\'e-Kuratowski limit of the $\Phi_p$-convex hull of two
points $x,y$ exists in $\mathbb R^n$ and we give an algebraic
characterization\footnote{ The Kuratowski-Painlev\'e lower limit of
the sequence of sets $\{A_n\}_{n\in\mathbb{N}}$, denoted
$Li_{n\to\infty}A_n$, is the set of points $p$ for which there
exists a sequence $\{p_n\}$ of points such that $p_n\in A_n$ for all
$n$ and $p = \lim_{n\to\infty}p_n$ ; a sequence $\{A_n\}_{n\in
\mathbb{N}}$ of subsets of $\Real^m$ is said to converge, in the
Kuratowski-Painlev\'e sense, to a set $A$ if $Ls_{n\to\infty}A_n = A
= Li_{n\to\infty}A_n$, in which case we write $A =
Lim_{n\to\infty}A_n$.}.

In this paper,  a weaker definition is proposed in line with the
 algebraic structure above introduced.

\begin{defn}\label{bconv} A subset $C$ of $\Real^n$ is $\mathbb B^\sharp$-convex if and only if for all
$t\in [0,1]$ and all $(x,y)\in C\times C$, $x\boxplus ty\in C$.
\end{defn}

Notice that for all $n$-dimensional orthant $K$ of $\mathbb R^n$, a
$\mathbb B$-convex subset of $K$ is $\mathbb B^\sharp$-convex. It is
shown that the following definitions of $\mathbb B^\sharp$-convexity
are equivalent.

\begin{prop} \label{comb} For all subset $C$ of $\mathbb R^n$, the following claims are equivalent:

\noindent $(a)$ $C$ is a $\mathbb B^\sharp$-convex subset of
$\Real^n$.
\medskip

\noindent $(b)$  For all $(x_1,...,x_m)\in C^m$ we have:
$$ \Big\{\bigboxplus_{i\in [m]}t_ix_i: \max_{i\in [m]} t_i=1,t\in [0,1]^m\Big\}\subset C.$$

\end{prop}
{\bf Proof:} Let us prove that $(a)$ implies $(b)$. If $C$ is
$\mathbb B^\sharp$-convex, this property is true for $m=2$. Suppose
it is true at rank $m$ and let us prove that it is true at rank
$m+1$. In other words, assume that for all $(x_1,...,x_m)\in C^m$ we
have: $ \Big\{\bigboxplus_{i\in [m]}t_ix_i: \max_{i\in [m]}
t_i=1,t\in [0,1]^m\Big\}\subset C$, we need to prove that if
$(x_1,...,x_m, x_{m+1})\in C^{m+1}$ then for all $t\in [0,1]^{m+1}$
such that $ \max_{i\in [m+1]} t_i=1$ we have $\bigboxplus_{i\in
[m+1]}t_ix_i\in C$. To establish this property, we use Proposition
\ref{preassoc} which implies that if $(x_1,...,x_m, x_{m+1})\in
C^{m+1}$ then, for all $t\in [0,1]^{m+1}$ such that $ \max_{i\in
[m+1]} t_i=1$, we have
$$\bigboxplus_{i\in [m+1]}t_ix_i=
\bigboxplus_{i\in [m+1]}\Big[t_ix_i\boxplus \Big(\bigboxplus_{j\in
[m+1]\backslash \{i\}}t_jx_j\Big)\Big].\quad (\star)$$
 For all $i$ set  $t_i^\star=\max\{t_j:j\in [m+1]\backslash \{i\}\}$.
 It follows that
 $$\bigboxplus_{i\in [m+1]}t_ix_i
 =\bigboxplus_{i\in [m+1]}\Big[t_ix_i\boxplus t_i^\star\Big(\bigboxplus_{j\in [m+1]\backslash \{i\}}
  (\frac{t_j}{t_i^\star})\,x_j\Big)\Big].$$
$(i)$ By definition, if $t_i^\star=\max\{t_j:j\in [m+1]\backslash
\{i\}\}<1$ then, since $ \max_{j\in [m+1]} t_j=1$, we have
$t_i=1$. Moreover, $\max_{j\in [m+1]\backslash
\{i\}}(\frac{t_j}{t_i^\star})=1$ and since, by hypothesis, the
property is assumed to be true at rank $m$, it follows that
  $\bigboxplus_{j\in [m+1]\backslash \{i\}}
  (\frac{t_j}{t_i^\star})\,x_j\in C.$ Hence, we deduce  that
  $$x_i\boxplus t_i^\star\Big(\bigboxplus_{j\in [m+1]\backslash \{i\}}
  (\frac{t_j}{t_i^\star})\,x_j\Big) \in C.$$
  $(ii)$ If $t_i^\star=1$ then there is some $i_0\in [m+1]\backslash\{i\}$ such that $t_{i_0}=1$ and, by hypothesis,
  it follows that $\bigboxplus_{j\in [m+1]\backslash \{i\}}t_jx_j \in C$. Furthermore, since $t_i\in [0,1]$
  we deduce from $(a)$
  that $t_ix_i\boxplus \Big(\bigboxplus_{i\in [m+1]\backslash\{i\}}t_ix_i\Big)\in C$.
   For all $i$, set $\Lambda_i= t_ix_i\boxplus \Big(\bigboxplus_{j\in [m+1]\backslash \{i\}}t_jx_j\Big)$.
   We have proven that, for each $i$,
  $\Lambda_i\in C$. Moreover, from Propositions \ref{preassoc} and  Lemma  \ref{assoc}, the $\Lambda_i$'s belong to a $n$-dimensional orthant
 $K$ and, then, can be composed associatively using the operation
 $\boxplus.$ Thus, we deduce from $(\star)$ that  $\bigboxplus_{i\in [m+1]}t_ix_i\in C$ which ends the
  proof of $(b)$. The converse inclusion is immediate. $\Box$\\

\begin{prop}\label{obvious} $(a)$ The emptyset, $\Real^n$, as well as all the singletons are
$\mathbb{B}^\sharp$-convex.
\medskip\noindent
$(b)$ If $\{D_\delta : \delta\in\Delta\}$ is an arbitrary family of
$\mathbb{B}^\sharp$-convex sets then $\bigcap_{\lambda}D_\delta$ is
$\mathbb{B}^\sharp$-convex.
\medskip\noindent
$(c)$ If $\{D_\lambda : \delta\in\Delta\}$ is a family of
$\mathbb{B}^\sharp$-convex sets such that $\forall \delta_1,
\delta_2\in\Delta$ \, $\exists \delta_3\in\Delta$ such that
$D_{\delta_{1}} \cup D_{\delta_{2}} \subset D_{\delta_{3}}$ then
$\bigcup_{\delta}D_{\delta}$ is  $\mathbb{B}^\sharp$-convex. $(d)$
If $C$ a $\mathbb{B}$-convex subset of $\Real_+^n$ then it is
$\mathbb{B}^\sharp$-convex.
\end{prop}

Given a set $S\subset\Real^n$ there is, according to $(a)$ above, a
$\mathbb{B}^\sharp$-convex set which contains $S$; by $(b)$ the
intersection of all such $\mathbb{B}^\sharp$-convex sets is
$\mathbb{B}$-convex; we call it the $\mathbb{B}^\sharp$-convex hull
of $S$ and we write $\mathbb{B}^\sharp[S]$ for the
$\mathbb{B}^\sharp$-convex hull of $S$.

\begin{prop}\label{closure} The following properties hold:

\par\noindent $(a)$ $\mathbb{B}^\sharp[\emptyset] = \emptyset$, $\mathbb B^\sharp[\Real^n] =
\Real^n$, for all $x\in\Real^n$,  $\mathbb B^\sharp[\{x\}] = \{x\}$.

\noindent$(b)$ For all $S\subset\Real^n$, $S\subset \mathbb
B^\sharp[S]$ and $\mathbb B^\sharp[ \, [\mathbb B^\sharp[S] \, ] =
\mathbb B^\sharp[S]$.

\noindent $(c)$ For all $S_1, S_2\subset\Real^n$, if $S_1\subset
S_2$ then $\mathbb B^\sharp[S_1] \subset \mathbb B^\sharp[S_2]$.

\noindent $(d)$ For all $S\subset\Real^n$, $\mathbb B^\sharp[S] =
\bigcup\{ \mathbb B^\sharp[A]: A \hbox{ is a finite subset of }
S\}$.

\noindent $(e)$  A subset $L\subset\Real^n$ is $\mathbb
B^\sharp$-convex if and only if, for all finite subset $A$ of $L$,
$\mathbb B^\sharp[A]\subset L$.
\end{prop}

\subsection{Intermediate Points and Copositivity}\label{fund}

A set of  points  we term the {\bf intermediate points} is
introduced to characterize  the $\mathbb B$-convex hull of two
points on the whole Euclidean vector space. For all $(x,y)\in
\Real^n\times \Real^n$, let us consider the map $\gamma:
\Real^n\times \Real^n\times [0,+\infty]\longrightarrow \Real$
defined by:
\begin{equation}\label{gamma}\gamma(x,y,t)=(\max\{1,t\})^{-1}\big(x
\boxplus {t}y\big),\quad \text{for all } t\geq 0
\end{equation} and by $\gamma(x,y,+\infty)=y$. For all $(x,y)\in
\Real^n\times \Real^n$, let  $ \mathcal I(x,y)$ be the subset
defined by $ \mathcal I(x,y)=\{i\in [n]: x_iy_i<0\}$ and let
$n(x,y)$ be its cardinal. Remark that $\gamma(x,y,0)=x$.

For all $i\in \mathcal I(x,y)$ and all $t_{i}^\star \in \Real_{++}$
a point $\gamma_{}\in \Real^n$ is called a {\bf $i$-intermediate
point} between $x$ and $y$ if there is some $t_{i}^\star\in
]0,+\infty[$ such that
\begin{equation}
\gamma_i (x,y,t_{i}^\star):=\big(\gamma (x,y,t_{i}^\star)\big)_i=0.
\end{equation}

\begin{lem}\label{gam}
Let $\gamma: \Real^n\times \Real^n\times
(\Real_{+}\cup\{+\infty\})\longrightarrow \Real^n$ be the map
defined in \eqref{gamma}. Suppose that $\mathcal
I(x,y)\not=\emptyset$. We have the the following properties:

\noindent $(a)$  For all $i\in \mathcal I(x,y)$, the map $t\mapsto
\gamma_i(x,y,t)$ has a uniqueness zero
$t_{i}^{\star}=-\frac{x_i}{y_i}=|\frac{x_i}{y_i}|>0$ and there is a
uniqueness $i$-intermediate point $$\gamma
(x,y,t_{i}^\star)=\Big(\frac{|y_i|}{\max\{|x_i|,|y_i|\}}\,x\Big)
\boxplus \Big(\frac{|x_i|}{\max\{|x_i|,|y_i|\}}\,y\Big). $$

\noindent $(b)$ For all $i\in \mathcal I(x,y)$ and all $t\geq 0$

$$\gamma_i(x,y,t)=\left\{\begin{matrix}\max\{1,t\}^{-1}x_i&\text{if}&t<-\frac{x_i}{y_i}\\
0&\text{if}&t=-\frac{x_i}{y_i}\\
t\max\{1,t\}^{-1}y_i&\text{if}&t>-\frac{x_i}{y_i}\end{matrix}\right.$$

\noindent $(c)$ If  $t>\max\{1, -\frac{x_i}{y_i}\}$ then
$\gamma_i(x,y,t)=y_i$. Moreover, if $t<\min\{1,-\frac{x_i}{y_i}\}$
then $\gamma_i(x,y,t)=x_i$.

\noindent $(d)$  $\lim_{t\longrightarrow 0}\gamma^{}(x,y,t)=x$ and
$\lim_{t\longrightarrow +\infty}\gamma^{}(x,y,t)=y$.

\end{lem}
{\bf Proof:} $(a)$ For all  $i\in \mathcal I(x,y)$, we have
$x_iy_i<0$, which implies that
$-\frac{x_i}{y_i}=|\frac{x_i}{y_i}|>0$. Moreover,
$\gamma_i(x,y,t)=0$ if and only if $ \Big(\max\{1,t\}^{-1}x_i\Big)
\boxplus \Big(\max\{1,t\}^{-1}ty_i\Big)=0$. Since this is
equivalent to $x_i\boxplus ty_i=0$, we deduce that
$t_i^\star=-\frac{x_i}{y_i}$ is the uniqueness positive zero of
the equation $\gamma_i(x,y,t)=0$. Moreover
$\gamma_i(x,y,0)=x_i\not=0$ and $\gamma_i(x,y,+\infty)=y_i\not=0$,
which ends the proof. $(b)$ The case $t=-\frac{x_i}{y_i}$ is an
immediate consequence of $(a)$. Assume that
$t\not=-\frac{x_i}{y_i}$. In such a case, one has
$$\gamma_i(x,y,t)=\left\{\begin{matrix}\max\{1,t\}^{-1}x_i&\text{if}&|\max\{1,t\}^{-1}x_i|>|t\max\{1,t\}^{-1}y_i|\\
t\max\{1,t\}^{-1}y_i&\text{if}&|\max\{1,t\}^{-1}x_i|<|t\max\{1,t\}^{-1}y_i|\end{matrix}\right..$$
Since $|\frac{x_i}{y_i}|=-\frac{x_i}{y_i}$, it follows that
$$\gamma_i(x,y,t)=\left\{\begin{matrix}\max\{1,t\}^{-1}x_i&\text{if}&t<-\frac{x_i}{y_i}\\
t\max\{1,t\}^{-1}y_i&\text{if}&t>-\frac{x_i}{y_i}\end{matrix}\right..$$

$(c)$ If $t>\max\{1, -\frac{x_i}{y_i}\}$, then, from $(b)$ one has
$\gamma_i(x,y,t)=t\max\{1,t\}^{-1}y_i$. Moreover $\max\{1,t\}=t$.
Therefore $\gamma_i(x,y,t)=y_i$. If $t<\min\{1,
-\frac{x_i}{y_i}\}$, then
 $\gamma_i(x,y,t)=\max\{1,t\}^{-1}x_i$. Moreover
$\max\{1,t\}=1$. Therefore $\gamma_i(x,y,t)=x_i$. $(d)$ Suppose
that $j\notin \mathcal I(x,y)$, then there is $\epsilon\in
\{-1,1\}$ such that $\Big(\max\big \{1,t\big \}^{-1}x\Big)_j
\boxplus \Big(t\max\big\{1,t \big\}^{-1}y\Big)_j\in
\epsilon\Real_+$. It follows that the map $t\mapsto
\gamma_j(x,y,t)$ is continuous. Hence, we clearly have
$\lim_{t\longrightarrow 0}\gamma_j^{}(x,y,t)=x_j$ and
$\lim_{t\longrightarrow +\infty}\gamma_j^{}(x,y,t)=y_j$. Suppose
now that $i\in \mathcal I(x,y)$. From $(c)$
$\lim_{t\longrightarrow 0}\gamma_i^{}(x,y,t)=x_i$
and $\lim_{t\longrightarrow +\infty}\gamma_i^{}(x,y,t)=y_i$ which ends the proof. $\Box$\\

Notice that it may happen that there are two indexes $i,k\in
\mathcal I(x,y)$ such that $\gamma (x,y,t_{i}^\star)=\gamma
(x,y,t_{k}^\star)$. Let $\Theta(x,y)=\{0,+\infty,
-\frac{x_i}{y_i}:i\in \mathcal I(x,y)\}$. If $\mathcal
I(x,y)=\emptyset$, then $\Theta(x,y)=\{x,y\}$.

\begin{expl}
Let  $x=(4,2)$, $x'=(3,4)$, $x''=(-\frac{7}{2}, 3)$ and
$y=(-2,-3)$ four points of $\mathbb R^2$. Clearly, we have
$\mathcal I(x,y)=\mathcal I(x',y)=\{1,2\}$. Let us denote
$\gamma_{i}$ the intermediate points between $x$ and $y$. We have
$\gamma_{1}=\max\{4,2\}^{-1}\Big(2(4,2)\boxplus
4(-2,-3)\Big)=(0,-3).$ The second intermediate point between $x$
and $y$ is
 $\gamma_{2}=\max\{2,3\}^{-1}\Big(3(4,2)\boxplus
2(-2,-3)\Big)=(4,0).$ Let us denote $\gamma_{i}'$ the intermediate
points between $x'$ and $y$. We have
$\gamma_{1}'=\max\{3,2\}^{-1}\Big(2(3,4)\boxplus
3(-2,-3)\Big)=(0,-3).$ Following a similar procedure, we get
$\gamma_{2}'=(\frac{9}{4},0)$. Finally, we have $\mathcal
I(x'',y)=\{2\}$, and using similar notations we obtain
$\gamma_{2}''=(-\frac{21}{8},0)$.
\end{expl}

\begin{lem} \label{intermediate}For all $x,y\in \mathbb R^n$ such that $\mathcal
I(x,y)\not=\emptyset$, let $\{i_m\}_{m\in [n(x,y)]}$ be an index
sequence of $\mathcal I(x,y)$ such that
$$ i_m\in \arg\max _{i\in \mathcal I(x,y)}\big \{-\frac{x_{i}}{y_{i}}: i\geq
m\big\}  \quad \text {for all} \quad m\in [n(x,y)].$$ We have for
all $ m\in [n(x,y) -1]$:
\begin{equation*}\label{conseq} \gamma
\big(x,y,-\frac{x_{i_{m}}}{y_{i_{m}}}\big)\boxdot \gamma
\big(x,y,-\frac{x_{i_{m+1}}}{y_{i_{m+1}}}\big)\geq 0.
\end{equation*} Moreover $ x\boxdot \gamma
\big(x,y,-\frac{x_{i_1}}{y_{i_1}}\big)\geq 0$, $ \gamma
\big(x,y,-\frac{x_{i_{n(x,y)}}}{y_{i_{n(x,y)}}}\big)\boxdot y\geq
0$.

\end{lem}
{\bf Proof:}  Suppose that $j\notin \mathcal I(x,y)$. In such a
case $x_jy_j\geq 0$. It follows that there is some $\epsilon\in
\{-1,1\}$ such that $x_j,y_j\in \epsilon \Real_+$. Therefore
$\gamma_j(x,y,t)=\big[\big(\max\big \{1,t\big \}^{-1}x\big)
\boxplus \big(t\max\big\{1,t \big\}^{-1}y\big)\big]_j\in \epsilon
\Real_+$, for all $t\geq 0$. Hence if $\gamma(x,y,t_ {i_m}^\star)$
and $\gamma(x,y,t_{i_{m+1}} ^\star)$ are two intermediate points,
 then for all $j\in [n]\backslash \mathcal I(x,y)$ one has
$$\big(\gamma_j(x,y,t_{i_{m}}^\star)\big)\big(\gamma_j(x,y,t_{i_{m+1}}^\star)\big)\geq
0.$$ Suppose now that $ i\in \mathcal I(x,y)$. Set
$t_{i}^\star=-\frac{x_{i_{}}}{y_{i_{}}}$ for all $i\in \mathcal
I(x,y)$. By construction $\{t_{i_m}^\star\}_{m\in [\,n(x,y)\,]}$
is a nondecreasing sequence of $\Real_{++}$. Since
$t_{i_{m}}^\star$ and $ t_{i_{m+1}}^\star$ are two consecutive
terms of this sequence, for all $i\in \mathcal I(x,y)$ we have
$t_{i}^\star\notin \;] t_{i_{m}}^\star, t_{i_{m+1}}^\star[$. Thus
one has either  $ t_i^\star \leq t_{i_{m}}^\star \leq
t_{i_{m+1}}^\star$ or  $ t_i^\star \geq t_{i_{m +1}}^\star \geq
t_{i_{m}}^\star$. From  Lemma \ref{gam}.b, {\small
\begin{align}&\gamma_i(x,y,t_{i_m}^\star)\;\gamma_i(x,y,t_{i_{m+1}}^\star)\nonumber \\&=\left\{\begin{matrix}
\max\{1,{t_{i_m}^\star}\}^{-1}\max\{1,t_{i_{m+1}}^\star\}^{-1}(x_i)^2\geq
0& \text{if}&t_i^\star>
t_{i_{m+1}}^\star\geq t_{i_m}^\star\\
0&\text{if}&t_i^\star\in \{t_{i_m}^\star,t_{i_{m+1}}^\star\}\\
{t_{i_m}^\star}{t_{i_{m+1}}^\star}\max\{1,{t_{i_m}^\star}\}^{-1}\max\{1,{t_{i_{m+1}}^\star}\}^{-1}(y_i)^2\geq
0&\text{if}&t_i^\star<
 t_{i_m}^\star\leq t_{i_{m+1}}^\star
\end{matrix}\right.\end{align}}
It follows that for all $m\in [n(x,y)-1]$, $\gamma
(x,y,t_{i_m})\boxdot \gamma (x,y,t_{i_{m+1}})\geq 0$. Let us prove
that $x\boxdot \gamma (x,y,t_{i_1})\geq 0$ and $ \gamma
(x,y,t_{i_{n(x,y)+1}})\boxdot y\geq 0$. Since for all $i\in
\mathcal I(x,y)$ $t_{i_{1}}^\star\leq t_{i}^\star\leq
t_{i_{n(x,y)}}^\star$ we have
$$x_i\;\gamma_i(x,y,t_{i_{1}}^\star)=\left\{\begin{matrix}\max\{1,t_{i_{1}}^\star\}^{-1}(x_i)^2\geq 0&
\text{if}&t_i^\star>
t_{i_{1}}^\star\\
0&\text{if}&t_i^\star=t_{i_{1}}^\star
\end{matrix}\right.$$
and
$$\gamma_i(x,y,t_{i_{n(x,y)}}^\star)\;y_i=\left\{\begin{matrix}
0&\text{if}&t_i^\star= t_{i_{n(x,y)}}^\star\\
{t_{i_{n(x,y)}}^\star}\max\{1,{t_{i_{n(x,y)}}^\star}\}^{-1}(y_i)^2\geq
0&\text{if}&t_i^\star< t_{i_{n(x,y)}}^\star
\end{matrix}\right.$$

Since $t_i^\star=-\frac{x_i}{y_i}>0$ for all $i\in \mathcal
I(x,y)$, this ends the proof.
$\Box$\\

Set $t_{i_0}^\star=0$, $t_{i_{n(x,y)+1}}^\star=+\infty$ and
$t_{i_m}^\star=-\frac{x_{i_{m}}}{y_{i_{m}}}$ for all $m\in
[n(x,y)]$. A  sequence $\{t_{i_m}^\star\}_{m=0}^{n(x,y)+1}$ of
$\Theta (x,y)$ satisfying the conditions of Lemma \ref{intermediate}
is called an {\bf intermediate sequence}.

One can then establish the following inclusion, whose the proof is
immediate.

\begin{prop}\label{BInclus} Let $C$ be a $\mathbb B^\sharp$-convex
set of $\Real^n$. Suppose that $x,y\in C$, and let
$\{t_{i_m}^\star\}_{m=0}^{n(x,y)+1}$ be an intermediate sequence
of $\Theta (x,y)$. Then, for all $m\in [n(x,y)]$,
$\gamma_{}(x,y,t_{i_m}^\star)\in C$. Furthermore
$$\bigcup_{m=0}^{n(x,y)}\mathbb
B\Big[\gamma_{}(x,y,t_{i_m}^\star),\gamma_{}(x,y,t_{i_{m+1}}^\star)\Big]\subset
C.$$

\end{prop}
{\bf Proof:} For all $t\geq 0$,
$\max\big\{\max\{1,t\}^{-1},t\max\{1,t\}^{-1}\}\big\}=1$.
Consequently, since $x,y\in C$, all the intermediate points lie in
$C$. Since for all $m\in \{0,1,...,n(x,y)\}$
$\gamma_{}(x,y,t_{i_m}^\star)$ and
$\gamma_{}(x,y,t_{i_{m+1}}^\star)$ are copositive, it follows that
$\mathbb
B\Big[\gamma_{}(x,y,t_{i_m}^\star),\gamma_{}(x,y,t_{i_{m+1}}^\star)\Big]\subset
C$ for all $m$, which ends the proof. $\Box$\\

The next results will be useful in the remainder of the paper.

\begin{lem}\label{identity}For all $(a,b,c,d)\in \Real^4$, if $\big(a\boxplus b\big)\big(c\boxplus d\big)\geq 0$ then
$$\big(a\boxplus b\big)\boxplus\big(c\boxplus d\big)=\digamma_{[4]}(a,b,c,d).$$\end{lem}
{\bf Proof:} We first assume that  $(a\boxplus b)(c\boxplus d)=0$.
In such a case, one has  either $a\boxplus b=0$ or $c\boxplus d=0$.
Suppose, for example, that $a\boxplus b=0$. Then $\big(a\boxplus
b\big)\boxplus\big(c\boxplus d\big)=0\boxplus\big(c\boxplus
d\big)=c\boxplus d.$ Moreover,
$\digamma_{[4]}(a,b,c,d)=\digamma_{[4]}(a,-a,c,d)=c\boxplus d$ and
the equality holds true. The proof is similar in the case where
$c\boxplus d=0$.

Suppose now that $(a\boxplus b)(c\boxplus d)>0$. Then $a\boxplus
b\not=0$ and $c\boxplus d\not=0$ and from Lemma \ref{symcopos}, we
deduce the result. $\Box$\\

For all $u,v,w,z\in \Real^n$, let us denote

\begin{equation}
u\boxplus v\boxplus w\boxplus
z=\big(\digamma_{[4]}(u_1,v_1,w_1,z_1),\cdots,\digamma_{[4]}(u_n,v_n,w_n,z_n)\big).
\end{equation}

\begin{prop}\label{BInclus}
For all $x,y\in \Real^n$, let $\{t_{i_m}^\star\}_{m=0}^{n
(x,y)+1}$ be an intermediate sequence of $\Theta (x,y)$. Then
\begin{align*}\bigcup_{m=0}^{n (x,y)}\mathbb
B&\Big[\gamma_{}(x,y,t_{i_m}^\star),\gamma_{}(x,y,t_{i_{m+1}}^\star)\Big]\\&\subset
\Big \{t x\boxplus r x\boxplus s y\boxplus w y:
\max\{t,r,s,w\}=1,t,r,s,w\geq 0\Big\}.\end{align*}

\end{prop}
{\bf Proof:} We have just to show that if  $\gamma,\gamma'\in
\Gamma(x,y)$ and $\gamma \boxdot \gamma'\geq 0$, then $\mathbb
B[\gamma,\gamma']\subset \Big \{t x\boxplus r x\boxplus s y\boxplus
w y: \max\{t,r,s,w\}=1,t,r,s,w\geq 0\Big\}.$ Suppose that $z\in
\mathbb B\big[\gamma,\gamma'\big]$. Hence, by hypothesis, there are
$\alpha,\alpha'\in \Real_+$ such that $\max\{\alpha,\alpha'\}=1$ and
$z=\alpha \gamma \boxplus \alpha' \gamma'$. Since $\gamma$ and
$\gamma'$ are two intermediate points there exists $s,t,s',t'\geq 0$
with $\max\{s,t\}=1$ and $\max\{s',t'\}=1$ and such that $\gamma=s
x\boxplus t y$ and $\gamma'=s' x\boxplus t' y$. It follows that
$z=[\alpha (s x\boxplus t y)]\boxplus [\alpha' (s'x\boxplus t'
y)]=[(\alpha s x)\boxplus (\alpha t y)]\boxplus [
(\alpha's'x)\boxplus (\alpha't' y)]$. Since $\gamma$ and $\gamma'$
are copositive, for all $i\in [n]$, $[ (\alpha s
 x_i)\boxplus (\alpha t y_i)][(\alpha' s'
 x_i)\boxplus (\alpha' t' y_i)]\geq 0$. We deduce from Lemma \ref{identity}  that
 $$\big[(\alpha s x_i)\boxplus (\alpha t y_i)\big]\boxplus \big[
(\alpha's'x_i)\boxplus (\alpha't' y_i)\big ]=( \alpha s x_i)\boxplus
(\alpha t y_i)\boxplus ( \alpha' s' x_i)\boxplus (\alpha' t'y_i).
$$
It follows that $z=( \alpha s) x\boxplus (\alpha t) y\boxplus (
\alpha' s') x\boxplus (\alpha' t') y$. Moreover, one has
$\max\{\alpha s, \alpha t, \alpha' s', \alpha't' \}=1$. Hence,
$z\in \Big \{t x\boxplus r x\boxplus s y\boxplus w y:
\max\{t,r,s,w\}=1,t,r,s,w\geq 0\Big\}$,
which ends the proof. $\Box$\\

\subsection{Some Topological Properties}

In the following we show that $\mathbb B^\sharp$-convex sets have a
path-connected structure. This we do using the intermediate function
and focusing on the copositive case.
\begin{prop}
Let $\gamma: \Real^n\times \Real^n\times
[0,+\infty]\longrightarrow \Real^n$ be the  map defined in
\eqref{gamma}. Suppose that $x$ and $y$ are copositive. Then:

\noindent $(a)$ The map  $t\mapsto \gamma(x,y,t)$ is continuous on
$\Real_{+}$.

\noindent $(b)$ One has  $\lim_{t\rightarrow
0}\gamma(x,y,t)=\gamma(x,y,0)=x$ and $\lim_{t\rightarrow
\infty}\gamma(x,y,t)=\gamma(x,y,+\infty)=y.$

\noindent $(c)$ For all copositive pairs $(x,y)$, we have $ \gamma
\big(x,y,[0,+\infty]\big)=\mathbb B[x,y]$.

\end{prop}
{\bf Proof:} $(a)$ The maps $t\mapsto \max\{1,t\}^{-1}$ is
continuous over $\Real_+$. Since $x$ and $y$ are copositive, for
all $i$ there is some $\epsilon_i\in \{-1,1\}$ such that
$$\gamma_i(x,y,t)=\epsilon_i\max\big\{\epsilon_i \max\{1,t\}^{-1}x_i,\epsilon_it\max\{1,t\}^{-1}y_i\big \}. $$
Consequently, each map $\gamma_i(x,y,\cdot)$ is continuous in $t$
and the result follows.

$(b)$ Clearly $\lim_{t\rightarrow 0}\max\{1,{t}\}^{-1}=1$ and
$\lim_{t\rightarrow 0}t\max\{1,t\}^{-1}=0$. Hence, we have
$\lim_{t\rightarrow 0}\gamma(x,y,t)=x$. Moreover,
$\lim_{t\rightarrow \infty}\max\{1,{t}\}^{-1}=0$ and
$\lim_{t\rightarrow \infty}t\max\{1,t\}^{-1}=1$. Consequently,
$\lim_{t\rightarrow \infty}\gamma(x,y,t)=y$.

 $(c)$ Since $x$ and $y$ are copositive, we
have $\mathbb B[x,y]=\{t x\boxplus s y: t,s\in [0,1], \max\{t,s\}=1
\}$. Since for all $t\geq 0$
$\max\big\{\max\{1,t\}^{-1},t\max\{1,t\}^{-1}\big\}=1$, we deduce
that $ \gamma (x,y,{\mathbb R}_{+})\subset \mathbb B[x,y]$. However,
from \cite{bh},  $\mathbb B[x,y]$ is a closed subset of $\Real^n$.
From $(a)$ and $(b)$, we deduce that $  \gamma
\big(x,y,[0,+\infty]\big)\subset \mathbb B[x,y]$. Let us show the
converse inclusion. By definition, we have $$\mathbb B[x,y]=\big\{t
x\boxplus  y: t\in [0,1]\big\}\cup \big\{x\boxplus s y: s \in
[0,1]\big\}.$$ Suppose that $0<t\leq 1$,  and set $t'=t^{-1}$. We
have
 $tx\boxplus y=\big(\max\{1,{{t'}\}^{-1}}\big)x\boxplus \big( t'\max\{1,{{t'}\}^{-1}}\big) y\in
 \gamma\big(x,y,[0,+\infty]\big)$. If  $t=0$ then
 $x=\gamma(x,y,0)$.
 Furthermore, if  $0\leq s\leq 1$, then
$x\boxplus s y=\big(\max\{1,{s }\}^{-1}x\big)\boxplus
\big(s\max\{1,{s}\}^{-1}y\big)\in \gamma(x,y,[0,+\infty])$, which
proves the converse inclusion.
 $\Box$\\

\begin{cor}\label{connect} Let $a$ and $b$ two real numbers with $a<b$. Let $h:[a,b]\longrightarrow \bar{\mathbb R}_+$
be an homeomorphism such that $h(a)=0$ and $h(b)=\infty$. Let
$\xi_h:\Real^n\times \Real^n\times [a,b]$ be the map defined by
$\xi_h(x,y,s)=\gamma(x,y,h(s))$. If $x$ and $y$ are copositive,
then the map $s\mapsto \xi_h(x,y,s) $ is continuous. Moreover
$\xi_h(x,y,a)=x$, $\xi_h(x,y,b)=y$ and $\xi_h\big(x,y,[a,b]\big)
=\mathbb B[x,y]$.

\end{cor}

It is shown below that a $\mathbb B^\sharp$-convex set is
path-connected.

\begin{prop} \label{precont} A non
empty $\mathbb B^\sharp$-convex  of $\Real_{}^{n}$ is
path-connected.
\end{prop}
{\bf Proof:} We first establish that for all $x,y\in \Real^n$, there
exists a continuous map $\xi (x,y,\cdot): [0,1]\longrightarrow
\Real^n$ such that   $$\xi (x,y,[0,1])=\bigcup_{m=0}^{n(x,y)}\mathbb
B\Big[\gamma_{}(x,y,t_{i_m}^\star),\gamma_{}(x,y,t_{i_{m+1}}^\star)\Big],$$
with $\xi(x,y,0)=x$ and $\xi(x,y,1)=y$. Let
$\{t_{i_m}^\star\}_{m=0}^{n(x,y)+1}$ be an intermediate sequence of
$\Theta(x,y)$  from $x$ to $y$. Set $x=\gamma (x,y,0)$ and $y=\gamma
(x,y,\infty)$.  From corollary \ref{connect}, for all $m\in
\{0,...,n(x,y)\}$ there is collection of continuous maps $$\eta_m:
[t_{i_m}^\star,t_{i_{m+1}}^\star]\longrightarrow \mathbb
B\Big[\gamma (x,y,t_{i_m}^\star),\gamma
(x,y,t_{i_{m+1}}^\star)\Big]$$ such that
$\eta_m([t_{i_m}^\star,t_{i_{m+1}}^\star])=\mathbb B\Big[\gamma
(x,y,t_{i_m}^\star),\gamma (x,y,t_{i_{m+1}}^\star)\Big]$,
$\eta(t_m)=\gamma (x,y,t_{i_m}^\star)$ and $\eta_m(t_{m+1})=\gamma
(x,y,t_{i_{m+1}})$. Let $\xi:\Real^n\times \Real^m\times
[0,1]\longrightarrow \Real^n$ be the map defined by
$\xi(x,y,t)=\eta_m(t)$ for all $t\in [t_m,t_{m+1}]$. Clearly, the
map $t\mapsto \xi(x,y,t)$ is continuous and we have $\xi(x,y,0)=x$
and $\xi(x,y,1)=y$. Moreover, from Proposition \ref{fond} one has
$\xi(x,y,[0,1])=\bigcup_{m=0}^{n(x,y)}\mathbb
B\Big[\gamma_{}(x,y,t_{i_m}^\star),\gamma_{}(x,y,t_{i_{m+1}}^\star)\Big]$.
Since   for all $x, y\in L$, $\xi(L\times
L\times [0, 1])\subset L$, this ends the proof. $\Box$\\

\subsection{Separation of Copositive $\mathbb B$-Convex Sets}

We say that two subsets  $C_1$ and $C_2$ of $\Real^n$ are
copositive if for all $(x_1,x_2)\in C_1\times C_2$, $x_1\boxdot
x_2\in \Real_+^n$. In this subsection, it is shown that the inner
product $(x,y)\mapsto \langle x,y\rangle_\infty=\bigboxplus_{i\in
[n]}
 x_iy_i$ can be used to separate two copositive $\mathbb B$-convex sets. For all $u,v\in \Real$, let us define the binary operation
\begin{equation*}\label{Bform}u\smile v=
\left\{\begin{matrix}
v &\hbox{ if } &|v|& > &|u|\\
\min \{u,v\}&\hbox{ if }&|u|&=&|v|\\
v& \hbox{ if }& |u|&<&|v|.\end{matrix}\right.\end{equation*} An
elementary calculus shows that $u\boxplus v=\frac{1}{2}\Big(u\smile
v-\big[(-u)\smile (-v)\big] \Big)$. It has been established in
\cite{bh3} that
 the set $\mathbb R$ equipped with the semilattice operation $\smile$ and the usual multiplication $\cdot$ by positive real
 numbers is a  semimodule over the semifield of positive real numbers $\mathbb R_+$.
 Furthermore, both $({\mathbb R}_+, \smile, \cdot)$ and $({\mathbb R}_-, \smile, \cdot)$ are sub-semimodules
 isomorphic to  $({\mathbb R}_+, \max, \cdot)$; the isomorphisms are, respectively, given by the inclusion,
 $u\mapsto u$, and the negative of the inclusion, $u\mapsto -u$.

  Given $m$ elements $u_1, \cdots, u_m$ of $\Real$, not all of which are $0$, let $I_+$, respectively $I_-$,
 be the set of indices for which $0 < u_i$, respectively $u_i < 0$. We can then write
 $u_1\smile\cdots\smile u_m = (\smile_{i\in I_+}u_i) \smile (\smile_{i\in I_-}u_i)
 = (\max_{i\in I_+}u_i)\smile (\min_{i\in I_-}u_i)$ from which we have

 \begin{equation}\label{manysmiles}
 u_1\smile\cdots\smile u_m =
\left\{
\begin{array}{lcc}
\max_{i\in I_+}u_i  &\hbox{if}   & I_- = \emptyset \hbox{ or } \max_{i\in I_-}\vert u_i\vert <  \max_{i\in I_+}u_i    \\
\min_{i\in I_-}u_i  &  \hbox{if} &  I_- = \emptyset \hbox{ or } \max_{i\in I_+}u_i <  \max_{i\in I_-}\vert u_i\vert \\
 \min_{i\in I_-}u_i  & \hbox{if}  &    \max_{i\in I_-}\vert u_i\vert  = \max_{i\in
 I_+}u_i.
\end{array}
\right.
 \end{equation}
We define a ${\mathbb B}$-form on ${\mathbb R}^n_+$ as a map $f:
{\mathbb R}^n_+\to\Real$ such that, for all
 $u_1, \cdots, u_m$ in $ {\mathbb R}^n_+$ and all $t_1, \cdots, t_m$ in $ {\mathbb R}_+$
$
 f(t_1u_1\vee\cdots \vee t_mu_m) = t_1f(u_1)\smile\cdots\smile t_mf(u_m).
$. It has been shown in \cite{bh3} that a map $f: {\mathbb
R}^n_+\to\Real$ is a ${\mathbb B}$-form if and only if there exists
$(a_1, \cdots, a_n)\in{\mathbb R}^n$, necessarily unique, such that,
for all $(x_1, \cdots, x_n)\in{\mathbb R}^n_+$,
 \begin{equation}\label{eqcarBform}
 f(x_1, \cdots, x_n) = a_1x_1\smile\cdots\smile a_nx_n.
 \end{equation}
Moreover for all ${\mathbb B}$-forms $f : {\mathbb R}_+^{n}\to\Real$
and all real numbers $c$ :
\begin{equation}\label{id1Bf}
\hbox{if } 0\leq c \hbox{ then } f(x)\leq c \hbox{ if and only if }
\max_{i\in I_+}\{a_ix\}\leq  \max_{i\in I_ -}\{-a_ix, c\}
\end{equation}
and
\begin{equation}\label{id2Bf}
\hbox{if } c\leq 0 \hbox{ then } f(x)\leq c \hbox{ if and only if }
\max_{i\in I_+}\{a_ix, -c\}\leq  \max_{i\in I_ -}\{-a_ix\}.
\end{equation}

For all $c\in \Real$ and all sunset $I$ of $[n]$ the map $y \mapsto
\max_{i\in I}\{y_i,c\} $ is continuous over $\Real^n$. Therefore for
all $c \in \Real$, $f^{-1}\left(\,]-\infty,c]\right)=\left\{x\in
\Real^n: f(x)\leq c\right\}$ is closed. It follows that  a ${\mathbb
B}$-form is lower semi-continuous.

The largest (smallest) lower (upper) semi-continuous minorant
(majorant) of a map $h$ is said to be the lower (upper)
semi-continuous regularization of $h$.

\begin{prop}\label{Reg}
\noindent Let $f$ be a $\mathbb B$-form defined by $ f(x_1, \cdots,
x_n) = a_1x_1\smile\cdots\smile a_nx_n, $ for some $a\in \Real^n$.
Then $f$ is the lower semi-continuous regularization of the map
$x\mapsto \langle a,x\rangle_{\infty}=\bigboxplus_{i\in [n]}a_ix_i$.
\end{prop}
\noindent{\bf Proof:} Suppose that $\varphi_a:\Real^n\longrightarrow
\Real$ is the lower semi-continuous regularization of $\langle
a,x\rangle_{\infty}$. First, remark that for all $x\in \Real^n$:
$$f(x)\leq \langle a,x\rangle_{\infty}.$$

Therefore, all we need to prove is that $\varphi_a(x)=f(x)$. By
definition, since $f$ is lower semi-continuous, we have for all
$x\in \Real^n$
$$f(x)\leq \varphi_a(x)\leq \langle a,x\rangle_{\infty}.$$
Let $I_-^a=\{i\in [n]: a_ix_i<0\}$. If $I_-^a=\emptyset$ then
$\langle a,x\rangle_{}^{\infty}=\max_{i=1\cdots n}|a_ix_i|=
\max_{i=1\cdots n}\{a_ix_i\}$. Moreover by definition $f(x)=
\max_{i=1\cdots n}\{a_ix_i\}=\langle a,x\rangle_{\infty}$.
Consequently, since $f(x)\leq \varphi_a(x)\leq \langle
a,x\rangle_{\infty}$, we deduce that $f(x)= \langle
a,x\rangle_{\infty}=\varphi_a(x).$

Suppose now that $I_-^a\not=\emptyset$ and pick some $i_0\in I_-^a$.
By hypothesis, we have $a_{i_0}\not=0$. Now, let
$\{x_k\}_{k\mathbb{N}}$ be the sequence defined as:
$$x_{k,i}=\left\{\begin{matrix}
x_i&\text{ if } i\not=i_0\\
x_{i_0}+\frac{1}{a_{i_0}k}&\text{ if } i=i_0.
\end{matrix}\right.$$
Hence, since $a_{i_0}x_{i_0}<0$ and $i_0\in I_-^a$ we have $\langle
a,x_k\rangle_{\infty}=a_{i_0}x_{i_0}-\frac{1}{k}=-\max_{i=1\cdots
n}|a_ix_i|-\frac{1}{k}$. Thus:
$$\lim_{k\longrightarrow
\infty}\langle a,x\rangle_{\infty}=\lim_{k\longrightarrow
\infty}\left(-\max_{i=1\cdots
n}|a_ix_i|-\frac{1}{k}\right)=-\max_{i=1\cdots n}|a_ix_i|.$$
Moreover, since $\varphi_a$ is lower semi-continuous and
$\lim_{k\longrightarrow \infty}x_k=x$:
$$\liminf_{k\longrightarrow \infty}\varphi_a(x_k)\geq \varphi_a(x)$$

By hypothesis $\varphi_a$ is  the lower semi-continuous
regularization of $\langle a,\cdot\rangle_{\infty}$, thus, by
definition, $\langle a,x_k\rangle_{\infty}\geq \varphi_a(x_k)$.
Therefore:
$$-\max_{i=1\cdots n}|a_ix_i|=\lim_{k\longrightarrow
\infty}\langle a,x_k\rangle_{\infty}\geq \liminf_{k\longrightarrow
\infty}\varphi_a(x_k)\geq \varphi_a(x)$$ Hence $\varphi_a(x)\leq
-\max_{i=1\cdots n}|a_ix_i|$. However, since $I_-^a\not=\emptyset$,
$f(x)=-\max_{i=1\cdots n}|a_ix_i|$, and we deduce that:
$$\varphi_a(x)\leq f(x).$$
But since $\varphi_a$ is the lower semi-continuous regularization of
the map $x\mapsto \langle a,x\rangle_{\infty}$ and $f(x)\leq \langle
a,x\rangle_{\infty}\;\forall x\in \Real^n$, we  also have:
$$\varphi_a(x)\geq f(x).$$
Consequently, $\varphi_a(x)= f(x)$ which ends the proof.
$\Box$\\

In \cite{bh3} it was established that if $C_1$ and $C_2$ are
nonproximate ${\mathbb B}$-convex subsets of ${\mathbb R}^n_+$ then
there exists a ${\mathbb B}$-form $f : {\mathbb R}_+^{n}\to\Real$
such that $
 \sup_{x\in C_1}f(x) < \inf_{x\in C_2}f(x).$ In the following, this
 result is extending to the inner product $(x,y)\mapsto \langle x,y\rangle_\infty=\bigboxplus_{i\in [n]}
 x_iy_i$.

  \begin{prop}  If $C_1$ and $C_2$ are nonproximate copositive ${\mathbb B}$-convex subsets of $\Real^n$ then there exists some  $a\in {\mathbb R}^n$ such that
 \begin{equation*}
 \sup_{x\in C_1}\langle a,x\rangle_\infty < \inf_{x\in C_2}\langle a,x\rangle_\infty.
 \end{equation*}
 \end{prop}
 {\bf Proof:} If $C_1$ and $C_2$ are copositive, then they belong to the same $n$-dimensional orthant $K$ that is homeomorphic to $\Real_+^n$ using a suitable linear homeomorphism. Therefore, for sake of simplicity, we shall assume
that $K=\Real_+^n$. From \cite{bh3}, there is some $a\in \Real^n$
such
 that the map $x\mapsto a_1 x_1\smile \cdots \smile a_nx_n$ separates
 $C_1$ and $C_2$. This implies that $ \inf_{x\in C_2}f(x)>\sup_{x\in
 C_1}f(x)$. Since $f$ is the lower semi-continuous regularization of $\langle
 a,\cdot\rangle_\infty$, it follows that $\langle
 a,x\rangle_\infty\geq f(x)$ for all $x\in C_2$. Therefore $$\inf_{x\in C_2}\langle a,x\rangle_\infty\geq \inf_{x\in
 C_2}f(x). \quad (\star)$$ Let us consider the map $g: \Real^n\longrightarrow
 \Real$ defined for all $x\in \Real^n$ by $g(x)=-f(-x)$. Since $x\mapsto
 -x$ is continuous, $g$ is upper semi-continuous. Moreover, for all $x\in \Real^n$, $f(-x)\leq \langle a,-x\rangle_\infty=-\langle a,x\rangle_\infty$
  implies that $g(x)\geq \langle a,x\rangle_\infty$. From equations \eqref{id1Bf} and   \eqref{id2Bf}, for all real numbers $c$ :
\begin{equation*}
\hbox{if }  c \leq 0\hbox{ then } g(x)\geq c \hbox{ if and only if }
\max_{i\in I_+}\{-a_ix\}\leq  \max_{i\in I_ -}\{a_ix, -c\}
\end{equation*}
and
\begin{equation*}
\hbox{if } c\geq 0 \hbox{ then } g(x)\geq c \hbox{ if and only if }
\max_{i\in I_+}\{-a_ix, c\}\leq  \max_{i\in I_ -}\{a_ix\}.
\end{equation*}
Hence $g(x)< c$ if and only if $f(x)< c$ and $\sup_{x\in C_1}g(x)<
\inf_{x\in C_2}f(x)\leq \inf_{x\in C_2}\langle a,x\rangle_\infty$.
Since  $\langle a
,x \rangle_\infty\leq g(x)$ for all $x\in C_1$, this ends the proof from $(\star)$. $\Box$\\

\section{ Relation to a Limit of Linear Convexities}

\subsection{Intermediate Points of Order $p$}

For all natural number $p$, let us consider the map $\gamma^{(p)}:
\Real^n\times \Real^n\times [0,+\infty]\longrightarrow \Real$
defined by:
\begin{equation}\label{gammap}\gamma^{(p)}(x,y,t)=\Big(\frac{1}{1\stackrel{p}{+}t}\Big) x\stackrel{p}{+}
\Big(\frac{t}{1\stackrel{p}{+}t}\Big) y ,\quad \text{for all } t\geq
0
\end{equation} and by  $\gamma^{(p)}(x,y,+\infty)=y$.

\begin{lem}\label{propP}For al $p\in \mathbb N$, the map defined in \eqref{gammap} satisfies the following
properties. For all $x,y\in \Real^n$:\\
$(a)$ The map $t\mapsto
\gamma^{(p)}(x,y,t)$ is continuous over $\Real_{+}$;\\
$(b)$ We have $\lim_{t\longrightarrow
0}\gamma^{(p)}(x,y,t)=\gamma^{(p)}(x,y,0)=x$ and
$\lim_{t\longrightarrow +\infty}\gamma^{(p)}(x,y,t)=\gamma^{(p)}(x,y,+\infty)=y$;\\
 $(c)$  We have $\gamma^{(p)}
(x,y,[0,+\infty])=Co^p(x,y)$.
\\
 $(d)$  For all $t\in [0,+\infty]$, we have $\lim_{p\longrightarrow +\infty}\gamma^{(p)}
(x,y,t)=\gamma^{} (x,y,t)$.\\
 $(e)$  For all $t\in [0,+\infty]$,  $\gamma^{(p)}
(x,y,t)$ and $\gamma^{} (x,y,t)$ are copositive.
\end{lem}
{\bf Proof:} $(a)$ The $\varphi_p$ generalized sum is continuous.
Moreover, for all $t\geq 0$ the map  $t\mapsto
(1\stackrel{p}{+}t)^{-1}$ is continuous and positive.  $(b)$ follows
from the continuity and using the fact that $\lim_{t\longrightarrow
+\infty}\frac{t}{1\stackrel{p}{+}t}=1$. $(c)$ Let us show that  for
all $t\geq 0$ $\gamma^{(p)}(x,y,t)\in Co^p(x,y)$. It is easy to see
that $\Big(\frac{1}{1\stackrel{p}{+}t}\Big) \stackrel{p}{+}
\Big(\frac{t}{1\stackrel{p}{+}t}\Big) =1 $ and  it follows  that
$\gamma^{(p)}(x,y,t)\in Co^p(x,y)$ for all $t\geq 0$.
 Since $\gamma^{(p)}(x,y,\infty)=y$, we deduce that $\gamma^{(p)}
(x,y,[0,+\infty])\subset Co^p(x,y)$. Conversely, suppose that $z\in
Co^p(x,y)$. By hypothesis there is some $\theta\in [0,1]$ such that
$z= \theta x \stackrel{p}{+} (1\stackrel{p}{-}\theta) y$. If
$\theta\in \{0,1\}$ then either $z=x$ or $z=y$. Suppose that
$\theta\in ]0,1[$, then setting $t=\frac{1\stackrel{p}{-}
\theta}{\theta}$, we obtain $\theta=\frac{1}{1\stackrel{p}{+}t}$ and
$1\stackrel{p}{-} \theta=\frac{t}{1\stackrel{p}{+}t}$. Consequently,
$z\in \gamma(x,y,\Real_{++})$. Therefore, the converse inclusion is
true and we deduce that $\gamma^{(p)} (x,y,[0,+\infty])=Co^p(x,y)$.
$(d)$ If either $t=0$ or $t=\infty$ then this property obviously
holds true. Suppose that $t\in \Real_{++}$. For all $j\in [n]$, we
have \begin{align*}\lim_{p\longrightarrow
\infty}\gamma_j(x,y,t)&=\lim_{p\longrightarrow
\infty}\Big(\frac{x_j}{1\stackrel{p}{+}t}\Big) \stackrel{p}{+}
\Big(\frac{ty_j}{1\stackrel{p}{+}t}\Big)= \lim_{p\longrightarrow
\infty}\Big(\frac{x_j^{2p+1}+(ty_j)^{2p+1}}{1+t^{2p+1}}\Big)^{\frac{1}{2p+1}}\\
&=\lim_{p\longrightarrow
\infty}\frac{(x_j^{2p+1}+(ty_j)^{2p+1})^{\frac{1}{2p+1}}}{(1+t^{2p+1})^{\frac{1}{2p+1}}}
=\frac{x_j\boxplus ty_j}{\max\{1,t\}}=\gamma_j(x,y,t).\end{align*}
$(e)$ From Lemma \ref{compar}, we have for all $j\in [n]$,
$\big(x_j^{2p+1}+(ty_j)^{2p+1}\big)^{\frac{1}{2p+1}}\geq 0$ if and
only if $x_j\boxplus t y_j\geq 0$. Using distributivity of scalar
multiplication, we deduce $(e)$. $\Box$\\

For all $i\in \mathcal I(x,y)$, we say that
$\gamma^{(p)}(x,y,t_{i}^\star)$ is {\bf a $i$-intermediate point
of order $p$ }between $x$ and $y$ if
$\gamma_i^{(p)}(x,y,t_{i}^{\star})=0$.

\begin{cor}\label{rstructure} Let $x,y\in \mathbb R^n$ and suppose that $\mathcal
I(x,y)\not=\emptyset$. Then, for all $i\in \mathcal I(x,y)$ one
has the following properties:

\noindent $(a)$ For all $p\in \mathbb N$, there is a uniqueness
$i$-intermediate point of order $p$
$$\gamma^{(p)}(x,y,t_{i}^\star)=\Big(\frac{|y_i|}
{|x_i|\stackrel{p}{+}|y_i|}\Big)x\stackrel{p}{+}
\Big(\frac{|x_i|}{|x_i|\stackrel{p}{+}|y_{i}|}\Big)y,$$ with
$t_i^\star=-\frac{x_i}{y_i}=|\frac{x_i}{y_i}|$.\\
\noindent $(b)$ For all $p\in \mathbb N$ and all $i\in \mathcal
I(x,y)$, $\gamma_{}^{(p)}\big (x,y,-\frac{x_i}{y_i}\big)$ is a
$i$-intermediate point of order $p$ if and only if $\gamma \big
(x,y,-\frac{x_i}{y_i}\big )$ is a $i$-intermediate point.\\
\noindent $(c)$ Let $\gamma (x,y,t_i^\star)$ be a $i$-intermediate
point and let $\{\gamma^{(p)} (x,y,t_i^\star)\}_{p\in \mathbb N}$ be
a sequence of $i$-intermediate points of order $p$. Then
$\lim_{p\longrightarrow \infty}\gamma^{(p)} (x,y,t_i^\star)
=\gamma^{} (x,y,t_i^\star)$.

\noindent $(d)$ If $\{t_{i_m}^\star\}_{m=0}^{n(x,y)+1}$ is an
intermediate sequence of $\Theta (x,y)$ satisfying the conditions of
Lemma \ref{intermediate} with $t_{i_0}^\star=0$,
$t_{i_{n(x,y)+1}}^\star=+\infty$ and
$t_{i_m}^\star=-\frac{x_{i_{m}}}{y_{i_{m}}}$ for all $m\in
[n(x,y)]$, then:
$$
Co^p(x,y)=\bigcup_{m=0}^{n(x,y)}Co^p\Big(\gamma^{(p)}
\big(x,y,t_{i_m}^\star\big),\gamma^{(p)}
\big(x,y,t_{i_{m+1}}^\star\big)\Big).
$$

\end{cor}
{\bf Proof:} $(a)$ $\gamma_i(x,y,t_i)=0$ if and only if
$\Big(\frac{1}{1\stackrel{p}{+}t}\Big)x_i\stackrel{p}{+}
\Big(\frac{t}{1\stackrel{p}{+}t}\Big) y_i=0$ which is equivalent to
$t_i=-\frac{x_i}{y_i}=|\frac{x_i}{y_i}|$. $(b)$ and $(c)$ are two
immediate consequences of $(a)$ and Proposition \ref{propP}.d. $(d)$
By definition, for all $m$ $\gamma^{(p)}
\big(x,y,t_{i_m}^\star\big),\gamma^{(p)}\in Co^p(x,y)$. Therefore $
Co^p(x,y)\supset \bigcup_{m=0}^{n(x,y)}Co^p\Big(\gamma^{(p)}
\big(x,y,t_{i_m}^\star\big),\gamma^{(p)}
\big(x,y,t_{i_{m+1}}^\star\big)\Big). $ Moreover, since
$x=\gamma^{(p)} \big(x,y,t_{i_0}^\star\big)$ and $y=\gamma^{(p)}
\big(x,y,t_{i_{n(x,y)+1}}^\star\big)$, the converse inclusion holds.
$\Box$\\

 \subsection{Painlev\'e-Kuratowski Limit}

 From Lemma \ref{intermediate} the $\Phi_p$-convex
hull of $x$ and $y$ is the finite union of the $\Phi_p$-convex hull
of two consecutive intermediate points of order $p$. The sequence
which these intermediate points of order $p$ are arranged is
identical to the copositive sequence of the intermediate points.

For future reference, we gather in the lemma below some elementary
facts, which are a slight extension of a result established in
\cite{bh}.
\begin{lem} \label{prel}Let $K$ be a $n$-dimensional orthant of $\mathbb
R^n$. Let $A=\{x_1,...,x_m\}$ be a finite subset of $K$. For all
natural number $p$ let $A^{(p)}=\{x_1^{(p)},...,x_m^{p)}\}$ be a
finite collection of $m$ vectors in $K$.

\noindent $(a)$ If there exists an an increasing sequence of natural
numbers $\{p_k\}_{k\in \mathbb{N}}$ such that for $i=1,...,m$
$\lim_{k\longrightarrow \infty}x_i^{(p_k)}=x_i$,
 then:\\ $$\lim_{k\longrightarrow
\infty}\stackrel{\varphi_{p_k}}{\sum_{i\in
[m]}}x_i^{(p_k)}=\bigboxplus_{i\in [m]}x_i.$$ $(b)$ If for
$i=1,...,m$ $\lim_{p\longrightarrow \infty}x_i^{(p)}=x_i$, then $$
Lim_{p\to\infty}Co^{p}(A^{(p)}) =\Big\{\bigboxplus_{i\in [m]}t_ix_i:
\max_i t_i=1, t_i\geq 0\Big\}=\mathbb B[A].$$
\end{lem}

\noindent {\bf Proof:} $(a)$  Let $\Psi_K:K\longrightarrow
\Real_+^n$ be the function characterizing the $n$-dimensional
orthant $K$. By definition for all $x\in \Real^n$ one has
$\Psi_K(x)=(\epsilon_1x_1,...,\epsilon_nx_n)$ where $\epsilon_j\in
\{-1,1\}$ for all $j\in [n]$. By definition one has for all $j$:
$$\stackrel{\varphi_{p_k}}{\sum_{i\in
[m]}}x_{i,j}^{(p_k)}=\epsilon_j \stackrel{\varphi_{p_k}}{\sum_{i\in
[m]}}|x_{i,j}^{(p_k)}|.$$ From the Lemma 2.0.1.b established in
\cite{bh}, we have $\lim_{k\longrightarrow
+\infty}\stackrel{\varphi_{p_k}}{\sum_{i\in
[n]}}|x_{i,j}^{(p_k)}|=\max_i|x_{i,j}|$. Consequently, for all $j\in
[n]$
$$\lim_{k\longrightarrow
+\infty}\stackrel{\varphi_{p_k}}{\sum_{i\in
[m]}}x_{i,j}^{(p_k)}=\epsilon_j\max_i|x_{i,j}|=\bigboxplus_{i\in
[m]}x_{i,j},$$ which ends the proof.

$(b)$ We first establish that $\mathbb B[A]=\left\{\bigboxplus_{i=1,
\cdots, m}t_i x_i : t_i\in [0,1], \max_{i \in[ m]}
t_i=1\right\}\subset Li_{p\longrightarrow \infty} Co^p(A)$. Let
$y=t_{1}x_{1}\boxplus\cdots\boxplus t_{m}x_{m}$ with
$t_{1},\cdots,t_{m}\in [0,1]$ and $\max_{i\in [m]}t_{i}=1$. Define
$y_{p}\in Co^{p}\left(A\right)$ by
$$y^{(p)}=\frac{1}{t_{1}\stackrel{p}{+}\cdots\stackrel{p}{+}t_{m}}
\left(t_{1}\stackrel{p}{.}x_{1}^{(p)}\stackrel{p}{+}\cdots
\stackrel{p}{+}\rho_{m}\stackrel{p}{.}x_{m}^{(p)}\right)$$ Since
$x_{1}^{(p)},\cdots,x_{m}^{(p)}\in K$ and $\lim_{p\longrightarrow
\infty}\left(t_{1}\stackrel{p}{+}\cdots\stackrel{p}{+}t_{m}\right)=
\max_{i\in[ m]}t_{i}=1$ we deduce from $(a)$ that
$$\lim_{p\to\infty}y^{(p)} = \lim_{p\longrightarrow
\infty}\left(t_{1}x_{1}^{(p)}\stackrel{p}{+}\cdots
\stackrel{p}{+}t_{m}x_{m}^{(p)}
\right)=t_{1}x_{1}\boxplus\cdots\boxplus t_{m}x_{m}=y.$$ This
completes the first part of the proof.
\medskip
Next, we establish that $Ls_{p\longrightarrow \infty} Co^p(A)\subset
\mathbb B[A]$. Take $z\in Ls_{p\longrightarrow \infty} Co^p(A)$;
there is an increasing sequence $\{p_k\}_{k\in \mathbb{N}}$ and a
sequence of points $\{z_k\}_{k\in \mathbb{N}}$ such that $z_k \in
Co^{p_k}(A^{(p_k)})$ and $\lim_{k\longrightarrow\infty} z_{k}=z$.
Each $z_k$ being in $Co^{p_k}(A^{(p_k)})$,  we can write
$$z_{k}=
t_{k,1}x_{1}^{(p_k)}\stackrel{p_k}{+}\cdots\stackrel{p_k}{+}
t_{k,m}x_{m}^{(p_k)}.$$ Since $t_k = (t_{k, 1}, ..., t_{k, m})\in
[0, 1]^m$ one can extract a subsequence
$(t_{k_l})_{l\in\mathbb{N}}$ converges to a point $t^* =
(t_{1}^{*}, ..., t_{m}^{*})\in [0, 1]^m$. It follows that for all
$i\in [m]$ $\lim_{l\longrightarrow
\infty}t_{k_l}x_i^{(p_{k_l})}=x_i $. Furthermore, from $(a)$,
$\lim_{l\to\infty}\left(\sum_{i=1}^{m}t_{k, i}^{2p_{k_l} + 1}
\right)^{1/(2p_{k_l} +1)} = \max_i{t_i^*}=1 $. It follows that for
all $i\in [m]$ $\lim_{l\longrightarrow
\infty}t_{k_l}x_i^{(p_{k_l})}=t_i^*x_i $. From $(a)$ we deduce
that
 $x =
\bigboxplus_{i=1}^{m}t_{i}^{*}x_i$ with $\max_{ i \in
[m]}\{t_{i}^{*}\} = 1$. The first and the second part of the proof
show that $$Ls_{p\longrightarrow \infty} Co^p(A^{(p)})\subset
\mathbb B[A]\subset Li_{p\longrightarrow \infty} Co^p(A^{(p)})$$ and
this completes the proof since we always have the inclusion
$Li_{p\longrightarrow \infty}
Co^p(A^{(p)})\subset Ls_{p\longrightarrow \infty} Co^p(A^{(p)})\quad\Box$\\

In the following, it is proven  that $Co^\infty
(\{x,y\})=Lim_{p\longrightarrow +\infty}Co^p(\{x,y\})$. This means
that, given two points in the whole Euclidean vector space, the
Painlev\'e-Kuratowski limit of their generalized convex hull exists.
Moreover it is established that it
 has an algebraic description.   For the sake of simplicity let
$Co^\infty (x,y)$ and $Co^p (x,y)$ denote these convex hulls for all
$p\in \mathbb N$. Let us consider $\ell$ sequences of subsets of
$\Real^n$ $\{A_m^{(p)}\}_{p\in \mathbb N}$, $m\in [\ell]$.  If there
exists a subset $A_m$ of $\Real^n$ such that $Lim_{p\longrightarrow
\infty}A_m^{(p)}=A_m$ for all $m\in [\ell]$, then it is easy to show
that:
\begin{equation} \label{limun}
Lim_{p\longrightarrow \infty}\Big(\bigcup_{m\in
[\ell]}A_m^{(p)}\Big)= \bigcup_{m\in [\ell]}A_m.
\end{equation}

\begin{prop}\label{fond}
For all $x,y\in \Real^n$, let $\{t_{i_m}^\star\}_{m=0}^{n(x,y)+1}$
be an intermediate sequence of $\Theta (x,y)$. Then
$$Co^{\infty}(x,y)=Lim_{p\longrightarrow +\infty}
Co^p(x,y)=\bigcup_{m=0}^{n(x,y)}\mathbb
B\Big[\gamma_{}(x,y,t_{i_m}^\star),\gamma_{}(x,y,t_{i_{m+1}}^\star)\Big].$$

\end{prop}
{\bf Proof:} From Lemma  \ref{propP}.d, we have for all $i\in
\mathcal I(x,y)$
$$\lim_{p\longrightarrow
+\infty}\gamma_{}^{(p)}\big(x,y,-\frac{x_i}{y_i}\big)=\gamma_{}^{}\big(x,y,-\frac{x_i}{y_i}\big).$$
Moreover, from Lemma  \ref{propP}.e, for all $p\in \mathbb N$ and
all $i\in  \mathcal I(x,y)$,
$\gamma_{}^{(p)}\big(x,y,-\frac{x_i}{y_i}\big)$ and
$\gamma_{}^{}\big(x,y,-\frac{x_i}{y_i}\big )$ are copositive. Recall
that two vectors are copositive if their components have the same
sign. From Proposition \ref{intermediate}, for all $m\in [n(x,y)]$
$\gamma (x,y,t_{i_m}^\star)$ and $\gamma (x,y,t_{i_{m+1}}^\star)$
are copositive. Hence it follows that for all $m$,
$\gamma^{(p)}(x,y,t_{i_m}^\star)$ and $\gamma^{(p)}
(x,y,t_{i_{m+1}}^\star)$ are copositive.

From Proposition \ref{prel}, we have for all $m$

{\small
\begin{align*}Lim_{p\longrightarrow +\infty}
Co^p\Big(\gamma^{(p)}(x,y,t_{i_{m}}^\star)
,\gamma^{(p)}(x,y,t_{i_{m+1}}^\star\Big)=\mathbb
B\Big[\gamma^{}(x,y,t_{i_{m}}^\star)
,\gamma^{}(x,y,t_{i_{m+1}}^\star)\Big].
\end{align*}}

 Moreover, we have from Corollary \ref{rstructure}.d {\small
\begin{align*} Co^p(x,y)
=
\bigcup_{m=0}^{n(x,y)}Co^p\Big(\gamma^{(p)}(x,y,t_{i_{m}}^\star),\gamma^{(p)}(x,y,t_{i_{m+1}}^\star)\Big).
\end{align*}}
Hence, from equation \eqref{limun}, and Corollary \ref{rstructure}.c,  the result follows.$\Box$\\

This property has an immediate consequence.

\begin{prop}A subset $C$ of $\mathbb R^n$ is $\mathbb
B^\sharp$-convex if and only if for all $x,y\in \mathbb R^n$
$Co^\infty(x,y)\subset C$.
\end{prop}
\noindent{\bf Proof:}  This  is an immediate consequence of
Propositions \ref{BInclus} and \ref{BInclus}. $\Box$\\

Notice that it is not clear from the definition of $\mathbb
B^\sharp$-convex sets that, for an arbitrary couple $(x,y)$ of
$\Real^n\times \Real^n$, $Co^\infty (x,y)$ is $\mathbb
B^\sharp$-convex.

\begin{cor}A $\mathbb B$-convex subset  of $\mathbb R^n$ is $\mathbb
B^\sharp$-convex.
\end{cor}

In \cite{bh} it was established that a $\mathbb B$-convex subset
of $\Real^n$ is connected. A stronger property is established
below. It is shown that $\mathbb B$-convex sets are
path-connected.

\begin{prop} A non empty $\mathbb{B}$-convex of
$\Real_{}^{n}$ is path-connected.
\end{prop}
\noindent{\bf Proof:}  By definition a subset $L$ of
 $\Real^{n}$ is $\mathbb B$-convex if for all
finite subset $A\subset L$  we have $Co^{\infty}(\{x, y\})\subset
L$. This implies that  $Co^{\infty}(\{x,
y\})\subset L$ for all $x, y\in L$, which yields the result from Proposition \ref{precont}. $\Box$\\

 In the next statement, an algebraic characterization of $Co^\infty(x,y)$.
 To prove this, we
use the fact that the convex hull is not modified whenever one
consider several occurrences of a given point. For example, for all
$p\in \mathbb N$, one can equivalently write:
\begin{align}
Co^p(x,y)&=\Big\{\alpha x\stackrel{p}{+} \beta y: \alpha
\stackrel{p}{+} \beta =1, \alpha,\beta \geq
0\Big\}\\&=\Big\{\alpha_1 x\stackrel{p}{+}\alpha_2 x\stackrel{p}{+}
\beta_1 y\stackrel{p}{+} \beta_2 y: \alpha_1 \stackrel{p}{+}\alpha_2
\stackrel{p}{+} \beta_1 \stackrel{p}{+} \beta_2 =1, \alpha_i,\beta_i
\geq 0\Big\}.\nonumber
\end{align}

\begin{lem} \label{LiInclus} For all $x,y\in \Real^n$,
 $$
Co^\infty(x,y)=\Big \{t x\boxplus r x\boxplus s y\boxplus w y:
\max\{t,r,s,w\}=1,t,r,s,w\geq 0\Big\}.$$
\end{lem}
{\bf Proof:} Suppose that $z\in \Big \{t x\boxplus r x\boxplus s
y\boxplus w y: \max\{t,r,s,w\}=1,t,r,s,w\geq 0\Big\}$. By definition
there exists $\alpha_1,\alpha_2,\beta_1,\beta_2\geq 0$ with
$\max\{\alpha_1,\alpha_2,\beta_1,\beta_2\}=1$ and such that
$$z=\alpha_1 x\boxplus \alpha_2 x \boxplus \beta_1 y\boxplus \beta_2 y.$$
Define
$$z^{(p)}=\frac{1}{\alpha_1
\stackrel{p}{+}\alpha_2 \stackrel{p}{+} \beta_1 \stackrel{p}{+}
\beta_2 }\big(\alpha_1 x\stackrel{p}{+}\alpha_2 x\stackrel{p}{+}
\beta_1 y\stackrel{p}{+} \beta_2 y\big).$$ By construction
$z^{(p)}\in Co^p[x,y]$. Taking the limit on both sides yields from
Proposition \ref{LimMaxFunc}: \begin{align*}\lim_{p\longrightarrow
\infty}z^{(p)}&= \frac{1}{\max\{\alpha_1,\alpha_2,\beta_1, \beta_2\}
}\big(\alpha_1 x\boxplus\alpha_2 x\boxplus \beta_1 y\boxplus \beta_2
y\big)\\&=\alpha_1 x\boxplus\alpha_2 x\boxplus \beta_1 y\boxplus
\beta_2 y=z.\end{align*} Consequently, $z\in Li_{p\longrightarrow
+\infty}Co^p(x,y)$.
 $\Box$


\bibliographystyle{amsplain}
\bibliography{xbib}

\begin{thebibliography}{99}


\bibitem{adil}{\sc Adilov, G.
 and A.M. Rubinov}, {$\mathbb B$-convex sets and functions}, {\it Numerical Functional Analysis and Optimization}, 27, (2006) pp. 237-257.

\bibitem{adilYe}{\sc Adilov, G. and I. Yesilce}, {$\mathbb B^{-1}$-convex sets and $\mathbb B^{-1}$-measurable maps},
{\it Numerical Functional Analysis and Optimization}, 33 (2012) pp.
131-141.

\bibitem{adilYeTi}{\sc Adilov, G. , I. Yesilce and G. Tinaztepe}, Separation of $\mathbb B^{-1}$-Convex Sets by
$\mathbb B^{-1}$-Measurable Maps, {\it Journal of Convex Analysis},
To apear.


\bibitem{avr1}{\sc Avriel, M.}, {$R$-convex Functions},  {\it Mathematical Programming}, 2, (1972) pp. 309-323.

\bibitem{avr2}{\sc Avriel, M.}, {\it Nonlinear Programming: Analysis and Methods}, Prentice Hall, New Jersey, 1976.

\bibitem{ben}{\sc Ben-Tal, A.}, { On Generalized Means and Generalized Convex Functions},  {\it Journal of Optimization Theory and Applications},
21,  (1977) pp. 1-13.

\bibitem{bh} {\sc Briec. W. and C.D. Horvath}, {\it  $\mathbb B$-convexity},  Optimization,
Vol. 53 (2), (2004) pp. 103-127 .

\bibitem{bh2} {\sc Briec, W. and C.D. Horvath}, {  Nash points, Ky Fan inequality and
equilibria of abstract economies in Max-Plus and B-convexity}, {\it
Journal of Mathematical Analysis and Applications}, 341 (2008), pp.
188–199.

\bibitem{bh3} {\sc Briec, W. and C.D. Horvath}, { On
the separation of convex sets in some idempotent semimodules}, {\it
Linear Algebra and its Applications}, 435 (2011) pp. 1542-1548.






\bibitem{ms}{\sc Maslov, V.P. and S.N. Samborski (eds)}, {\it Idempotent Analysis},
Advances in Soviet Mathematics, American Mathematical Society,
Providence, 1992.

\bibitem{pat}{\sc Patriche, M. } { Fixed point and equilibrium theorems in a generalized
convexity framework}, {\it Journal of Optimization Theory and
Applications} 156 (2013), pp. 701-715.



\bibitem{rubi}{\sc Rubinov, A.}, {\it Abstract Convexity and Global Optimization}, Kluwer, 2000.



\end{thebibliography}

\end{document}